\documentclass[11pt,letterpaper,reqno]{amsart}
\usepackage{amsmath,amsthm,amsfonts,amssymb,mathtools,algpseudocode}
\usepackage{comment,bm}

\usepackage{bookmark,hyperref}
\hypersetup{pdfstartview={FitH}}

\newtheorem{Theorem}{{\bf Theorem}}[section]
\newtheorem{Example}[Theorem]{{\bf Example}}

\newtheorem{Proposition}[Theorem]{{\bf Proposition}}
\newtheorem{Definition}[Theorem]{{\bf Definition}}
\newtheorem{Lemma}[Theorem]{{\bf Lemma}}

\numberwithin{equation}{section}

\newcommand{\diag}{\mbox{\rm diag\ }}

\newcommand{\Ol}{\mbox{\Large $\mathcal O$}}

\begin{document}

\title[Jacobi Method for Symmetric Matrices of order $4$ under Parallel Strategies]{On the Global Convergence of the Jacobi Method for Symmetric Matrices of order $4$ under Parallel Strategies}

\author{Erna Begovi\'{c}~Kova\v{c}}\thanks{Erna Begovi\'{c}~Kova\v{c}, Faculty of Chemical Engineering and Technology, University of Zagreb, Maruli\'{c}ev trg 19, 10000 Zagreb, Croatia}
\author{Vjeran Hari}\thanks{Vjeran Hari, Department of Mathematics, Faculty of Science, University of Zagreb, Bijeni\v{c}ka 30, 10000 Zagreb, Croatia}
\thanks{This work has been fully supported by Croatian Science Foundation under the project 3670.}
\date{5 March 2017}


\begin{abstract}
The paper analyzes special cyclic Jacobi methods for symmetric matrices of order $4$. Only those cyclic pivot strategies that enable full
parallelization of the method are considered. These strategies, unlike the serial pivot strategies, can force the method to be very slow or
very fast within one cycle, depending on the underlying matrix. Hence, for the global convergence proof one has to consider two or three
adjacent cycles. It is proved that for any symmetric matrix $A$ of order~$4$ the inequality $S(A^{[2]})\leq(1-10^{-5})S(A)$ holds, where $A^{[2]}$ results from
$A$ by applying two cycles of a particular parallel method. Here $S(A)$ stands for the Frobenius norm of the strictly upper-triangular part of
$A$. The result holds for two special parallel strategies and implies the global convergence of the method under all possible fully parallel
strategies. It is also proved that for every $\epsilon>0$ and $n\geq4$ there exist a symmetric matrix $A(\epsilon)$ of order $n$ and a
cyclic strategy, such that upon completion of the first cycle of the appropriate Jacobi method the inequality $S(A^{[1]})> (1-\epsilon)S(A(\epsilon))$ holds.
\end{abstract}

\maketitle

\section{Introduction}

The Jacobi method applies a sequence of similarity transformations by plane rotations to a symmetric matrix in order to diagonalize it. The method can be described as an iterative process of the for
$$A^{(k+1)} = R_{k}^{T}A^{(k)} R_{k}, \quad k\geq0; \qquad\ A^{(0)}=A,$$
where $R_{k}$ are plane rotations and $A$ is a symmetric matrix of order $n$. The method is \textit{globally convergent} if, for each starting $A$, the generated sequence $(A^{(k)})$ converges to a diagonal matrix. Its global (asymptotic) convergence has been considered in \cite{for+hen-60,han-63,BP,hen+zim-68,SS89} (\cite{wil-62,har-91}) and its accuracy in \cite{dem+ves-92,drm+ves-04a,drm+ves-04b,mat-09}. A~one-sided version of the method has been studied in \cite{hes-58,ves+har-89} and the block versions in \cite{drm-07,H15,B}. There are many papers on Jacobi methods, and further references can be found within the bibliographies of the papers cited above.

At the step $k$ the method annihilates two off-diagonal elements of $A^{(k)}$, $a_{i(k)j(k)}^{(k)}$ and $a_{j(k)i(k)}^{(k)}$, $i(k)<j(k)$. The element $a_{i(k)j(k)}^{(k)}$ is the \emph{pivot element} while $i=i(k)$ and $j=j(k)$ are \emph{pivot indices}. The way of selecting the pivot pair at each step is called \emph{pivot strategy}.
The elements of $R_{k}$ are the same as in the identity matrix $I_{n}$, except for the elements at positions $(i,i)$, $(i,j)$, $(j,i)$, $(j,j)$, which are $\cos\varphi^{(k)}$, $-\sin\varphi^{(k)}$, $\sin\varphi^{(k)}$, $\cos\varphi^{(k)}$, respectively. The \emph{rotation angle} is determined by the known formula
\begin{equation}\label{kut1}
\tan2\varphi^{(k)} = \frac{2a_{ij}^{(k)}}{a_{ii}^{(k)}-a_{jj}^{(k)}}, \quad \varphi^{(k)}\in[-\pi/4,\pi/4],
\end{equation}
which implies
\begin{align}
a_{ii}^{(k+1)} &= a_{ii}^{(k)} +\tan\varphi^{(k)} a_{ij}^{(k)}, \label{aii'} \\
a_{jj}^{(k+1)} &= a_{jj}^{(k)} -\tan\varphi^{(k)} a_{ij}^{(k)}, \label{ajj'}
\end{align}
and
$$S^2(A^{(k+1)}) = S^2(A^{(k)})-(a_{ij}^{(k)})^2.$$
Here $S(X)$ stands for the \textit{off-norm} of a symmetric matrix $X$ of order $n$,
$$S(X) = \frac{\sqrt{2}}{2}\|X-\textrm{diag}(X)\|_F = \sqrt{\sum_{i=1}^{n-1}\sum_{j=i+1}^n x_{ij}^2}, \quad X=X^T = (x_{ij}).$$
In the definition (\ref{kut1}) of the rotation angle, we assume that $\varphi^{(k)}=0$ if $a_{ij}^{(k)}=0$ and $a_{ii}^{(k)}=a_{jj}^{(k)}$. It is the most natural assumption which can be rephrased as: if the pivot element is zero, just skip it.

Since the diagonal elements converge if the rotation angle is chosen as in the relation (\ref{kut1}) (see \cite{mas-95}), it is easy to show that the obtained sequence $(A^{(k)})$ converges to some diagonal matrix if and only if
\begin{eqnarray}\label{conv}
\lim_{k\rightarrow\infty}S(A^{(k)})=0.
\end{eqnarray}
Therefore, the method is globally convergent if (\ref{conv}) holds for any initial $A$. Since the sequence $(S(A^{(k)}))$ is nonincreasing, for the global convergence of the method it is sufficient to show that for any symmetric matrix $A$ we have
\begin{equation}\label{conv2}
S^{2}(A^{(\tau N)})\leq\gamma S^{2}(A), \quad0\leq\gamma<1, \ \tau\in\{1,2,3\}, \ N=\frac{n(n-1)}{2},
\end{equation}
where $\gamma$ and $\tau$ do not depend on $A$. Here we prove that the relation (\ref{conv2}) holds with $\tau=2$ or $\tau=3$, for the
case $n=4$ and for those cyclic strategies which enable parallel processing. For these strategies one cycle (or sweep) consists of three ``parallel steps''.

Why would one consider the Jacobi method for symmetric matrices of order $4$ when that problem can be solved directly? Jacobi method is known for its high relative accuracy on well behaved symmetric matrices, for its efficiency on nearly diagonal matrices and for its suitability for parallel processing. So, the natural choice of a pivot strategy for matrices of order $4$ is a parallel strategy. We have discovered that parallel strategies are very special. Depending on the underlying matrix, the reduction of the quantity $S(A)$ per sweep
can be extremely slow or fast. This knowledge can be used to improve the implementation of the algorithm. Finally, the Jacobi method for large symmetric positive definite matrices is nowadays implemented as one sided block algorithm. At each step the block algorithm has to solve the same eigenvalue problem but for much smaller matrix, typically of order $16$--$256$. For this purpose one can use an element-wise Jacobi method or one can accelerate it by using the block algorithm which solves a $4$ by $4$ eigenvalue problem at each step.

There are several comments related to the inequality (\ref{conv2}) and its proof. First, the proof presented here reveals that the reduction of the quantity $S(A)$ during one cycle can be arbitrary small. It sheds light to convergence failure of the cyclic Jacobi method discussed in \cite{BP}. We show that for every $\epsilon>0$ there is a starting matrix $A(\epsilon)$ and a cyclic Jacobi method such that upon completion of the first cycle the inequality $S(A^{[1]})>(1-\epsilon)S(A(\epsilon))$ holds. This fact is first proved for $n=4$ and then for any $n\geq4$. Hence the global convergence consideration for the general cyclic Jacobi method should scrutinize more than one cycle of the process.
Second, the presented result covers the most difficult part in the proof that every cyclic Jacobi method for symmetric matrices of order $4$ is globally convergent \cite{B,HB15a}.

The paper is divided into five sections and three appendices. In Section~\ref{sec2} we introduce notation and the basic concepts of the theory of equivalent strategies. We also recall some known convergence results. In Section~\ref{sec3} we concentrate on parallel strategies and introduce an auxiliary tool, a linear operator $\mathcal{T}_{A}$, which simplifies
the convergence analysis. The convergence result is formulated and proved for some trivial cases. Section~\ref{proof} is devoted to the global convergence proof and Section~\ref{sec5} to the construction of the above mentioned matrix $A(\epsilon)$ and to the proofs of the related results. Since the proofs of the main results are pretty complicated, we have moved all lengthy and technical proofs to appendices A, B and C. They are related to the results from sections~\ref{sec3}, \ref{proof} and \ref{sec5}, respectively.

Some of the results presented here can be found in the unpublished thesis \cite{B}.

\section{Basic concepts and notation}\label{sec2}

For the Jacobi method for symmetric matrices of order $n$, the \textit{pivot strategy} can be defined as a function $I:\mathbb{N}_{0}\rightarrow\mathbf{P}_{n}$, where $\mathbb{N}_{0}=\{0,1,2,3,\ldots\}$ and $\mathbf{P}_{n}=\big\{(i,j) \ \big| \ 1\leq i<j\leq n\big\}$. We say that at step $k$, $I$ selects the pivot pair $I(k)=(i(k),j(k))$ which lies in~$\mathbf{P}_{n}$. Let $I$ be a pivot strategy. If there is a positive integer $T$ such that $I(k+T)=I(k)$ for all $k\geq0$, we say that $I$ is \textit{periodic} with period $T$. If $T=N\equiv\frac{n(n-1)}{2}$ and $\{I(k) \ \big| \ 0\leq k\leq T-1\}=\mathbf{P}_{n}$, the pivot strategy is \textit{cyclic}.

For $S\subseteq\mathbf{P}_{n}$, let $\mathcal{\Ol}(S)$ denote the set of all finite sequences made of the elements of $S$, assuming that each pair from $S$ appears at least once
in each sequence from $\mathcal{\Ol}(S)$. Let $\mathcal{O}$ be a sequence of pairs from $\mathcal{\Ol}(S)$. An \textit{admissible transposition} on $\mathcal{O}$ is any transposition of two adjacent pairs from $\mathcal{O}$,
$$(i_r,j_r),(i_{r+1},j_{r+1}) \rightarrow (i_{r+1},j_{r+1}),(i_r,j_r),$$
provided that $\{i_{r},j_{r}\}\cap\{i_{r+1},j_{r+1}\}=\emptyset$. For such pairs we say that they \textit{commute}, or that they are \emph{disjoint}.
Two sequences $\mathcal{O}, \mathcal{O}' \in\mathcal{\Ol}(S)$ are called
\begin{itemize}
\item[(i)] \textit{Equivalent} if one can be obtained from the other by a finite number of admissible transpositions. Then we write $\mathcal{O}\sim\mathcal{O}'$.
\item[(ii)] \textit{Shift-equivalent} if $\mathcal{O}=[\mathcal{O}_{1},\mathcal{O}_{2}]$ and $\mathcal{O}'=[\mathcal{O}_{2},\mathcal{O}_{1}]$, where $[\mathcal{O}_{1},\mathcal{O}_{2}]$ stands for the concatenation of the sequences $\mathcal{O}_{1}$ and $\mathcal{O}_{2}$. We write $\mathcal{O}\stackrel{s}{\sim}\mathcal{O}'$.
\item[(iii)] \textit{Weakly equivalent} if one can find $\mathcal{O}_{1},\ldots,\mathcal{O}_{r-1}$ from $\mathcal{\Ol}(S)$ such that in the sequence $\mathcal{O}\equiv\mathcal{O}_{0},\mathcal{O}_{1},\ldots,\mathcal{O}_{r}\equiv\mathcal{O}'$, each pair of adjacent terms
    $\mathcal{O}_{i}$, $\mathcal{O}_{i+1}$, $0\leq i\leq r-1$, consists of either equivalent or shift-equivalent terms. In such a case, we write \mbox{$\mathcal{O}\stackrel{w}{\sim}\mathcal{O}'$}.
\end{itemize}
One can check that $\sim$, $\stackrel{s}{\sim}$ and $\stackrel{w}{\sim}$ are equivalence relations on $\mathcal{\Ol}(S)$. In our application we shall have $S=\mathbf{P}_{n}$.

Once these equivalence relations are defined on $\mathcal{\Ol}(\mathbf{P}_{n})$, they can easily be transferred to the set of cyclic pivot strategies. Here is the procedure.

Let $I$ be a cyclic pivot strategy. By $\mathcal{O} _{I} $ we mean the sequence of pairs $I (0)$, $I (1)$, $\:\ldots\:$, $I (N-1)$. Conversely, for $\mathcal{O}\in\mathcal{\Ol}(\mathbf{P}_{n})$, $\mathcal{O} = (i_{0},j_{0}),(i_{1},j_{1}),\ldots,(i_{N-1},j_{N-1})$, the cyclic strategy generated by $\mathcal{O}$ is defined by
$I_{\mathcal{O}}(k)=(i_{\omega(k)},j_{\omega(k)})$, provided that $k\equiv \omega(k) \ (\mathrm{mod}\ N)$, $0\leq\omega(k)\leq N-1$, $k\geq0$.
In other words, $I_{\mathcal{O}}(k)$ runs through $\mathcal{O}$ in the cyclic way as $k$ increases.

Two cyclic strategies $I$ and $I'$ are equivalent (we write $I\sim I'$), shift-equivalent ($I\stackrel{s}{\sim} I'$) and weakly equivalent ($I \stackrel{w}{\sim} I'$) if the same is true for the corresponding sequences $\mathcal{O}_{I}$ and~$\mathcal{O} _{I'}$. Note that for the shift-equivalent strategies we have $I'(k)=I(k+\sigma)$, $k\geq0$, for some shift $\sigma$, $0\leq\sigma\leq N-1$. (We can confine to nonnegative shifts since $I(k-\sigma)=I(k+N-\sigma)$.)

The importance of weakly equivalent cyclic strategies comes from the following result.

\begin{Theorem}\label{SSekv}\emph{\cite{SS89}}
If the Jacobi method converges for some cyclic strategy $I$, then it also converges for all strategies that are weakly equivalent to $I$.
\end{Theorem}

Note that Theorem~\ref{SSekv} also covers the cases of equivalent and shift-equivalent strategies. Another important result regarding the convergence under two weakly equivalent strategies is proved in \cite[Lemma~4.8]{H15}.

A cyclic strategy $I$ can be represented by the matrix $M_{I}=(m_{ij})$, where
$$m_{ij}=m_{ji}=k, \quad \text{if} \ I(k)=(i,j), \ i<j,$$
and $m_{ss}=-1$, $1\leq s\leq n$. Instead of $-1$, we shall display $*$ to indicate that the diagonal positions are not part of the pivot sequence (see (\ref{4x4})).
If $I=I_{\mathcal{O}}$, we shall also write $M_{\mathcal{O}}$.

\section{Parallel strategies in the case $n=4$}\label{sec3}

Let $A$ be a symmetric matrix of order $4$. Since the length of each $\mathcal{O}\in\mathcal{\Ol}(\mathbf{P}_{4})$ equals $\frac{4\cdot3}{2}=6$, each cyclic Jacobi method applies six steps within one cycle. Among all cyclic strategies a distinguished role is played by the ``parallel'' ones. They enable parallel processing, so the corresponding method will be called parallel Jacobi method (cf.~\cite{sam-71}). Each parallel Jacobi method for symmetric matrices of order $4$ applies three parallel steps within each cycle. Every parallel step consists of two consecutive steps which can be performed concurrently. This way, instead of six sequential steps, using a parallel pivot strategy, we apply three parallel steps within one cycle.

As we shall see, it will be sufficient to study just two cyclic pivot strategies $I_{1}$ and $I_{2}$, which have the following two-dimensional representations
\begin{equation}\label{4x4}
M_{I_{1}}=\left[\begin{array}{cccc}
* & 4 & 0 & 2 \\
4 & * & 3 & 1 \\
0 & 3 & * & 5 \\
2 & 1 & 5 & * \\
\end{array}\right] \quad \text{and} \quad
M_{I_{2}}=\left[\begin{array}{cccc}
* & 4 & 2 & 0 \\
4 & * & 1 & 3 \\
2 & 1 & * & 5 \\
0 & 3 & 5 & * \\
\end{array}\right],
\end{equation}
respectively. All cyclic strategies that can be fully parallelized, are shift equivalent to $I_{1}$ or $I_{2}$. Therefore, the convergence results for all parallel strategies follow
from the results for the strategies $I_{1}$ and $I_{2}$.

Consider the sets of pairs $\{(1,3),(2,4)\}$, $\{(1,4),(2,3)\}$, $\{(1,2),(3,4)\}$. Note that the pairs within braces commute. These are the only sets that contain commuting pairs and only they can define parallel Jacobi steps. From the first braces we see that the corresponding plane rotations $R(1,3,\varphi_{13})$ and $R(2,4,\varphi_{24})$ commute, and also their entries can be computed independently of each other. So, the corresponding Jacobi steps can be applied in parallel: first apply concurrently the left transformations and then
the right ones, or vice versa. This corresponds to the one parallel step which consists of two subsequent ordinary Jacobi steps. The same can be said for the steps corresponding to the other two braces. This leads us to parallel strategies, which we represent by the matrices
$$\mathbf{M}_1=\left[
    \begin{array}{cccc}
      * & 2 & 0 & 1 \\
      2 & * & 1 & 0 \\
      0 & 1 & * & 2 \\
      1 & 0 & 2 & * \\
    \end{array}
  \right] \quad \text{and} \quad
\mathbf{M}_2=\left[
    \begin{array}{cccc}
      * & 2 & 1 & 0 \\
      2 & * & 0 & 1 \\
      1 & 0 & * & 2 \\
      0 & 1 & 2 & * \\
    \end{array}
  \right], \ \text{respectively}.$$
Here, the matrix entries count the parallel steps and mark the pivot positions associated with them.

By inspecting all commuting pairs, we conclude that there are exactly six parallel strategies and they can be grouped into two clusters which are actually equivalent classes for the relation $\stackrel{s}{\sim}$. They are defined by the following orderings from $\mathcal{\Ol}(\mathbf{P}_{4})$, where $\mathcal{O}_{j}\stackrel{s}{\sim}\mathcal{O}_{j}'\stackrel{s}{\sim}\mathcal{O}_{j}''$, $j=1,2$,
\begin{footnotesize}
\begin{align*}
\mathcal{O}_{1} &= (1,3),(2,4),(1,4),(2,3),(1,2),(3,4), \quad
&\mathcal{O}_{2} = (1,4),(2,3),(1,3),(2,4),(1,2),(3,4), \\
\mathcal{O}_{1}' &= (1,2),(3,4),(1,3),(2,4),(1,4),(2,3), \quad
&\mathcal{O}_{2}' = (1,2),(3,4),(1,4),(2,3),(1,3),(2,4), \\
\mathcal{O}_{1}'' &= (1,4),(2,3),(1,2),(3,4),(1,3),(2,4), \quad
&\mathcal{O}_{2}'' = (1,3),(2,4),(1,2),(3,4),(1,4),(2,3).
\end{align*}
\end{footnotesize}
In order to prove the global convergence of the Jacobi method under all six parallel strategies, it is sufficient to prove it for the strategies $I_{1}$ and $I_{2}$. This follows from Theorem~\ref{SSekv}. Next, we show that the strategies $I_{1}$ and $I_{2}$ are closely connected, so that the method converges under one of them if and only if it converges under the other one. To this end, note that the matrices $\mathbf{M}_{1}$ and $\mathbf{M}_{2}$ are permutationally similar,
\begin{equation}\label{pequiv}
M_{2}=P^{T}M_{1}P,
\end{equation}
where $P=P_{12}$ or $P=P_{34}$. Here $P_{ij}$ is the transposition which interchanges rows (columns) $i$ and $j$ if a matrix is premultiplied (postmultiplied) by it. If (\ref{pequiv}) holds, we say that $I_{2}$ and $I_{1}$ are \emph{permutationally equivalent} (see~\cite{B,HB15a}).

\begin{Proposition}\label{tm:prva}
Let $A=(a_{ij})$ be a symmetric matrix of order $4$. Let $A^{(0)}=A,A^{(1)},\ldots$ be obtained by applying the cyclic Jacobi method defined by the strategy $I_{2}$ on $A$. Let $P=P_{12}$ or $P=P_{34}$, and let $\mathsf{A}^{(0)}=P^{T} AP,\mathsf{A}^{(1)},\ldots$ be obtained by applying the cyclic Jacobi method defined by the strategy $I_{1}$ on $\mathsf{A}^{(0)}$. Then $\mathsf{A}^{(2r)}=P^{T}A^{(2r)}P$, $r\geq0$.
\end{Proposition}

\begin{proof} The proof has been moved to~\ref{appA}.
\end{proof}

Thus, Proposition~\ref{tm:prva} implies that the Jacobi method converges under the strategy $I_{2}$ if and only if it converges under the strategy $I_{1}$. In particular, if the relation (\ref{conv2}) holds for the method defined by $I_{1}$, with some $\tau$ and $\gamma$, it holds for the method defined by $I_{2}$ with the same $\tau$ and
$\gamma$, and vice versa. Theorem~\ref{tm:main} below, shows that the relation (\ref{conv2}) holds for the strategy $I_{1}$ with $\tau= 2$ and $\gamma=1-10^{-5}$.

What can be said for the method under the strategies $I_{\mathcal{O}_{1}'}$, $I_{\mathcal{O}_{1}''}$ and $I_{\mathcal{O}_{2}'}$, $I_{\mathcal{O}_{2}''}$? For these strategies, the relation (\ref{conv2}) holds with the same $\gamma$ and with $\tau$ larger for $1$. In particular, for $\tau= 3$ and $\gamma=1-10^{-5}$. We shall show it for $I_{\mathcal{O}_{1}''}$. For the other three strategies the proof is similar.

Let us apply the Jacobi method defined by the strategy $I_{\mathcal{O}_{1}''}$ to a symmetric matrix~$A$ of order $4$, thus generating the sequence of matrices $A^{(0)}=A$, $A^{(1)}$, $A^{(2)},\ldots{}$. Let us consider $\tau+1=3$ cycles of the method. We display each second iterate, i.e. the iterates obtained after each of the first nine
parallel steps:
\begin{small}
\begin{align*}
A && \stackrel{(1,4),(2,3)}{\longrightarrow} && A^{(2)} && \stackrel{(1,2),(3,4)}{\longrightarrow} && A^{(4)} && \stackrel{(1,3),(2,4)}{\longrightarrow} && A^{(6)} \\
&& \stackrel{(1,4),(2,3)}{\longrightarrow} && A^{(8)} && \stackrel{(1,2),(3,4)}{\longrightarrow} && A^{(10)} && \stackrel{(1,3),(2,4)}{\longrightarrow} && A^{(12)} \\
&& \stackrel{(1,4),(2,3)}{\longrightarrow} && A^{(14)} && \stackrel{(1,2),(3,4)}{\longrightarrow} && A^{(16)} && \stackrel{(1,3),(2,4)}{\longrightarrow} && A^{(18)}.
\end{align*}\end{small}
We concentrate on the matrix $A^{(4)}$. If another Jacobi method is applied to $A^{(4)}$, the one defined by the strategy $I_{1}$, one obtains (after each two steps) the same matrices $A^{(6)},A^{(8)},A^{(10)},\ldots{}$. After two sweeps, one obtains the matrix $A^{(16)}$ and, if the relation (\ref{conv2}) holds for $I_{1}$ with $\tau=2$ and $\gamma<1$, then one has $S(A^{(16)}) \leq\gamma S(A^{(4)})$. Therefore, one obtains
$$S(A^{[3]}) =S(A^{(3 N)}) = S(A^{(18)})\leq S(A^{(16)}) \leq\gamma S(A^{(4)})\leq S(A),$$
proving the claim.

\subsection{The cyclic strategy $I_{1}$}

We focus on strategy $I_{1}=I_{\mathcal{O}_{1}}$ where
$$\mathcal{O}_{1} = (1,3),(2,4),(1,4),(2,3),(1,2),(3,4).$$
By this strategy, at the beginning of each cycle (except for the first cycle), the elements at the positions $(1,2)$ and $(3,4)$ are zero. Since we consider the global convergence, we can assume that the initial matrix already has the form
\begin{equation}\label{4x4matrica}
A = \left[\begin{array}{cccc}
a_{11} & 0 & a_{13} & a_{14} \\
0 & a_{22} & a_{23} & a_{24} \\
a_{13} & a_{23} & a_{33} & 0 \\
a_{14} & a_{24} & 0 & a_{44} \\
\end{array}\right].
\end{equation}

Let $[e_{1}\ e_{2}\ e_{3}\ e_{4}]$ denote the column partition of the identity matrix, and let
\begin{equation}\label{Q}
Q=[e_{1}\ e_{3}\ e_{4}\ {-}e_{2}]=\left[\begin{array}{rrrr}
1 & 0 & 0 & 0 \\
0 & 0 & 0 & -1 \\
0 & 1 & 0 & 0 \\
0 & 0 & 1 & 0 \\
\end{array}\right].
\end{equation}
The similarity transformation with $Q$ and with $Q^{T}$ has the following effect on the elements of a square matrix $X=(x_{rs})$,
\begin{equation}\label{Q'XQ}
Q^TXQ ={\small\left[
\begin{array}{rrrr}
      x_{11} & x_{13} & x_{14} & -x_{12} \\
      x_{31} & x_{33} & x_{34} & -x_{32} \\
      x_{41} & x_{43} & x_{44} & -x_{42} \\
      -x_{21} & -x_{23} & -x_{24} & x_{22} \\
    \end{array} \right]}, \qquad
QXQ^T = {\small\left[
\begin{array}{rrrr}
      x_{11} & -x_{14} & x_{12} & x_{13} \\
      -x_{41} & x_{44} & -x_{42} & -x_{43} \\
      x_{21} & -x_{24} & x_{22} & x_{23} \\
      x_{31} & -x_{34} & x_{32} & x_{33} \\
    \end{array} \right]}.
\end{equation}
Thus, for each $A$ we have $S(Q^{T}AQ)=S(A)$. We see a favorable movement of the elements lying at the pivot positions for the parallel steps. We can use it to define a new iterative process, closely related to the original Jacobi process, where the pivot elements always remain at the same positions. This will simplify the analysis.

Therefore, we introduce a linear operator which is comprised of the transformation corresponding to the first parallel step under $I_{1}$ followed by the similarity transformation with $Q$.

\begin{Definition}\label{T}
Let $\mathcal{S}_4$ denote the vector space of real $4$ by $4$ symmetric matrices. For $A\in\mathcal{S}_4$ let
$$\mathcal{T}_{A}(H)=(R(1,3,\phi)R(2,4,\psi)Q)^T HR(1,3,\phi)R(2,4,\psi)Q, \quad\ H\in \mathcal{S}_4,$$
where $R(1,3,\phi)$ and $R(2,4,\psi)$ are Jacobi rotations which annihilate the elements $a_{13}$ and $a_{24}$ of $A$, respectively, and $Q$ is defined by the relation (\ref{Q}). The rotation angles $\phi$, $\psi$ are from the interval $[-\frac{\pi}{4},\frac{\pi}{4}]$, so that the formulas (\ref{kut1})--(\ref{ajj'}) hold.

For $A\in\mathcal{S}_4$ let $\mathcal{T}(A)=\mathcal{T}_{A}(A)$ and for any $k\geq 0$
\begin{eqnarray*}
\mathcal{T}^k (A ) &=& \underbrace{\mathcal{T}(\mathcal{T}(\ldots(\mathcal{T}}(A)\ldots))), \quad \mathcal{T}^0(A)=A. \\
&&\hspace{6ex}k
\end{eqnarray*}
\end{Definition}

Thus, $\mathcal{T}_{A}: \mathcal{S}_{4}\mapsto\mathcal{S}_{4}$ is a linear operator. Note that if $a_{13}=0$ and $a_{24}=0$, then $\mathcal{T}_{A}$ reduces to the similarity transformation with the similarity matrix $Q$. The function $\mathcal{T}$ is not linear. However, it satisfies
$$S(\mathcal{T}^{k+1}(A)) \leq S(\mathcal{T}^{k}(A)), \quad k\geq0.$$
If $A\in\mathcal{S}_{4}$ is as in the relation (\ref{4x4matrica}) and $A'=\mathcal{T}(A)$, then we have
$$\mathcal{T}:\left[
    \begin{array}{cccc}
      a_{11} & 0 & a_{13} & a_{14} \\
      0 & a_{22} & a_{23} & a_{24} \\
      a_{13} & a_{23} & a_{33} & 0 \\
      a_{14} & a_{24} & 0 & a_{44} \\
    \end{array}
  \right]\mapsto\left[
    \begin{array}{cccc}
      a_{11}' & 0 & a_{13}' & a_{14}' \\
      0 & a_{22}' & a_{23}' & a_{24}' \\
      a_{13}' & a_{23}' & a_{33}' & 0 \\
      a_{14}' & a_{24}' & 0 & a_{44}' \\
    \end{array}
  \right],$$
with
\begin{subequations}\label{4_elementi}
\begin{align}
a_{11}' & =a_{11}+a_{13}\tan\phi, & \qquad a_{13}' &= a_{14}\cos\phi\cos\psi-a_{23}\sin\phi\sin\psi, \\
a_{22}' & = a_{33}-a_{13}\tan\phi, & \qquad a_{14}' & = -a_{14}\cos\phi\sin\psi-a_{23}\sin\phi\cos\psi, \\
a_{33}' & = a_{44}-a_{24}\tan\psi, & \qquad a_{23}' & =  -a_{14}\sin\phi\cos\psi-a_{23}\cos\phi\sin\psi, \\
a_{44}' & = a_{22}+a_{24}\tan\psi, & \qquad a_{24}' & = a_{14}\sin\phi\sin\psi-a_{23}\cos\phi\cos\psi.
\end{align}
\end{subequations}
The rotation angles $\phi\in[-\frac{\pi}{4},\frac{\pi}{4}]$ and $\psi\in[-\frac{\pi}{4},\frac{\pi}{4}]$ are determined by
\begin{equation}\label{4_kutevi}
\tan(2\phi)=\frac{2a_{13}}{a_{11}-a_{33}}, \qquad\tan(2\psi)=\frac{2a_{24}}{a_{22}-a_{44}}.
\end{equation}

First, we show that the repeated application of $\mathcal{T}$ to $A$ yields the matrices which are closely related to Jacobi iterations under the parallel strategy $I_{1}$.

\begin{Proposition}\label{tm:prop1}
Let $A\in\mathcal{S}_{4}$ and let $A^{(2k)}$ be obtained by applying $2k$ steps of the Jacobi method under the strategy $I_{1}$ to $A$. Then
\begin{equation}\label{Tprop}
\mathcal{T}^{k}(A)=(Q^{k})^{T} A^{(2k)}Q^{k}, \quad k\geq0.
\end{equation}
\end{Proposition}

\begin{proof} The proof is lengthy and technical, so we have moved it to~\ref{appA}.
\end{proof}

In particular, the relation (\ref{Tprop}) implies
\begin{equation}\label{STA}
S(\mathcal{T}^{k}(A))=S(A^{(2k)}), \quad k\geq0.
\end{equation}
We use Proposition~\ref{tm:prop1} to simplify the proof of the main result which follows.

\begin{Theorem}\label{tm:main}
Let $A\in\mathcal{S}_{4}$ be such that $a_{12}=0$, $a_{34}=0$ and let $A^{(12)}$ be obtained by applying $12$ steps of the Jacobi method under the strategy $I_{1}$ to $A$. Then
\begin{equation}\label{tmT}
S(A^{(12)})\leq(1-\epsilon)S(A),
\end{equation}
with $\epsilon= 10^{-5}$.
\end{Theorem}

Note that $12$ steps correspond to two sweeps of the method. Theorem~\ref{tm:main} ensures the global convergence of the method since the sequence of iterates $(S(A^{(l)}),l\geq0)$ is nonincreasing and its subsequence $(S(A^{(12t)}),t\geq0)$ converges to zero.

The proof of the main theorem is lengthy, hence we will devote the entire Section~\ref{proof} to it. However, we first provide a lemma that covers the special cases when more than two off-diagonal elements are equal to zero. Then the relation (\ref{tmT}) holds with much larger $\epsilon$ ($\epsilon=1$ or $1/2$).

\begin{Lemma}\label{4x4lema}
Let $A=(a_{rs})\in\mathcal{S}_{4}$ be such that $a_{12}=0$ and $a_{34}=0$. If
\begin{itemize}
\item[(i)] $a_{14}=0$ and $a_{23}=0$, then $A^{(2)}$ is diagonal.
\item[(ii)] $a_{13}=0$ and $a_{24}=0$, then $A^{(4)}$ is diagonal.
\item[(iii)] $a_{13}=0$, then  $S^2(A^{(4)})\leq\frac{1}{2}S^2(A)$.
\item[(iv)] $a_{24}=0$, then $S^2(A^{(4)})\leq\frac{1}{2}S^2(A)$.
\item[(v)] $a_{14}=0$, then $S^2(A^{(4)})\leq\frac{1}{2}S^2(A)$.
\item[(vi)] $a_{23}=0$, then $S^2(A^{(4)})\leq\frac{1}{2}S^2(A)$.
\end{itemize}
\end{Lemma}

\begin{proof} The proof has been moved to~\ref{appA}.
\end{proof}

\section{Proof of Theorem~\ref{tm:main}}\label{proof}

Let $\epsilon=10^{-5}$. Then the assertion (\ref{tmT}) of Theorem~\ref{tm:main} can be expressed in the form
\begin{equation}\label{tmT1}
S(\mathcal{T}^{6}(A))\leq(1-\epsilon)S(A).
\end{equation}
Instead of working with matrices $A^{(l)}$, $0\leq l\leq12$, we shall work with $B^{(k)}=\mathcal{T}^{k} (A)=(b_{rs}^{(k)})$, $0\leq k\leq6$.
Let $B\equiv B^{(0)}$, so that $B=A$ holds. If $S(B)=0$, then Theorem~\ref{tm:main} holds. We assume $S(B)>0$.

Contrary to the assertion of the theorem suppose that
\begin{equation}\label{4t}
S(B^{(6)}) > (1-\epsilon)S(B).
\end{equation}
We shall show that the relation (\ref{4t}) leads to a contradiction.

From Lemma~\ref{4x4lema}, we conclude that all off-diagonal elements of $B^{(k)}$ except for $b_{12}^{(k)}$ and $b_{34}^{(k)}$ are non-zero for $0\leq k\leq4$. Furthermore, by Lemma~\ref{4x4lema}(i), we have $|b_{14}^{(5)}|+|b_{23}^{(5)}|>0$, and since $S(B^{(k)})\leq S(B^{(k-1)})$, $k\geq1$, we have
\begin{equation}\label{4x4off}
S(B^{(k)}) > (1-\epsilon)S(B), \quad \text{for} \ 0\leq k\leq6.
\end{equation}
Let
\begin{equation}\label{4x4delta}
\delta_{k}=\frac{\sqrt{(b_{13}^{(k)})^{2}+(b_{24}^{(k)})^{2}}}{S(B)}, \quad k\geq0; \quad \delta=\delta_{0}.
\end{equation}
From our assumptions it follows that $\delta_{k}>0$ for at least $0\leq k\leq4$. Note that
\begin{equation}\label{new1}
S^{2}(B^{(k+1)}) = S^{2}(B^{(k)}) - (b_{13}^{(k)})^{2}-(b_{24}^{(k)})^{2} = (b_{14}^{(k)})^{2}+(b_{23}^{(k)})^{2}, \quad k\geq0.
\end{equation}
This implies
$$0< \delta_0^2+\delta_1^2+\delta_2^2+\delta_3^2+\delta_4^2+\delta_5^2 = 1-\frac{S^2(B^{(6)})}{S^2(B)} < 2\epsilon-\epsilon^2,$$
in particular
\begin{equation}\label{new3}
\delta_{k} < \sqrt{2\epsilon-\epsilon^{2}} < 0.0045, \quad 0\leq k\leq5.
\end{equation}
Formulas (\ref{4_elementi}) describe the transition from $B^{(k-1)}$ to $B^{(k)}$ for any $k\geq1$. In this transition we denote the angles $\phi$ and $\psi$ by $\phi_{k}$ and $\psi_{k}$, respectively. If we set $c_{\phi_{k}}=\cos\phi_{k}$, $c_{\psi_{k}}=\cos\psi_{k}$, $s_{\phi_{k}}=\sin\phi_{k}$, $s_{\psi_{k}}=\sin\psi_{k}$, $k\geq1$, then from the formulas (\ref{4_elementi}) we have
\begin{align*}
& \quad\big(b_{13}^{(k)}\big)^{2}+\big(b_{24}^{(k)}\big)^{2} = \\
& = \left( (b_{14}^{(k-1)})^{2}+(b_{23}^{(k-1)})^{2}\right) (c_{\phi_{k}}^{2} c_{\psi_{k}}^{2}+s_{\phi_{k}}^{2}s_{\psi_{k}}^{2}) - 4b_{14}^{(k-1)} b_{23}^{(k-1)}c_{\phi_{k}}c_{\psi_{k}}s_{\phi_{k}}s_{\psi_{k}} \\
& \geq\left( (b_{14}^{(k-1)})^{2}+(b_{23}^{(k-1)})^{2}\right)(c_{\phi_{k}}^{2} c_{\psi_{k}}^{2}+s_{\phi_{k}}^{2}s_{\psi_{k}}^{2}) - 4|b_{14}^{(k-1)}b_{23}^{(k-1)}c_{\phi_{k}}c_{\psi_{k}}s_{\phi_{k}}s_{\psi_{k}}| \\
& = \left( |b_{14}^{(k-1)}|-|b_{23}^{(k-1)}|\right) ^{2}(c_{\phi_{k}}^{2} c_{\psi_{k}}^{2}+s_{\phi_{k}}^{2}s_{\psi_{k}}^{2}) +
2|b_{14}^{(k-1)}b_{23}^{(k-1)}|\big(c_{\phi_{k}}c_{\psi_{k}}-|s_{\phi_{k}}s_{\psi_{k}}|\big)^{2} \\
& \geq\frac{1}{2}\left( |b_{14}^{(k-1)}|-|b_{23}^{(k-1)}|\right) ^{2}.
\end{align*}
Here, $c_{\phi_{k}}^{2} c_{\psi_{k}}^{2}+s_{\phi_{k}}^{2}s_{\psi_{k}}^{2}$ has been bounded by $\frac{1}{2}$, as in the proof of Lemma~\ref{4x4lema}(v). This implies
\begin{equation}\label{b14-b23}
\left| \,|b_{14}^{(k-1)}|-|b_{23}^{(k-1)}|\,\right| \leq \sqrt{2}\sqrt{(b_{13}^{(k)})^{2}+(b_{24}^{(k)})^{2}}=\sqrt{2}\delta_{k}S(B), \quad k\geq1.
\end{equation}
Together with (\ref{new3}), the relation (\ref{b14-b23}) yields
\begin{equation}\label{4x4a14-a23}
\big||b_{14}^{(k)}|-|b_{23}^{(k)}|\big| \leq \sqrt{2}\delta_{k+1}S(B) < 2\sqrt{\epsilon}S(B), \quad 0\leq k\leq4.
\end{equation}

\begin{Lemma}\label{tm:lema3}
Exactly one of the following two assertions holds:
\begin{itemize}
\item[(a)] $|b_{14}^{(k)}+b_{23}^{(k)}|\leq\sqrt{2}\delta_{k+1}S(B) < 2\sqrt{\epsilon}S(B)$, \quad $0\leq k\leq 4$,
\item[(b)] $|b_{14}^{(k)}-b_{23}^{(k)}|\leq\sqrt{2}\delta_{k+1}S(B) < 2\sqrt{\epsilon}S(B)$, \quad $0\leq k\leq 4$.
\end{itemize}
\end{Lemma}

\begin{proof}
Suppose that the two inequalities in (a) hold for some $k$, $0\leq k\leq4$. Then, because of the relations (\ref{new1}) and (\ref{4x4off}), we have
\begin{equation}\label{proofLemma3}
\frac{|b_{14}^{(k)}+b_{23}^{(k)}|}{|b_{14}^{(k)}|+|b_{23}^{(k)}|} < \frac{2\sqrt{\epsilon}S(B)}{\sqrt{(b_{14}^{(k)})^2+(b_{23}^{(k)})^2}}
= \frac{2\sqrt{\epsilon}S(B)}{S(B^{(k+1)})}
< \frac{2\sqrt{\epsilon}S(B)}{(1-\epsilon )S(B)}
= 2 \frac{\sqrt{\epsilon}}{1-\epsilon}<1,
\end{equation}
and therefore
\begin{align}
\nonumber|b_{14}^{(k)}-b_{23}^{(k)}| & =|b_{14}^{(k)}|+|b_{23}^{(k)}| \geq \sqrt{|b_{14}^{(k)}|^{2}+|b_{23}^{(k)}|^{2}}=S(B^{(k+1)}) \\
\label{proofLemma3a} & > (1-\epsilon)S(B)>2\sqrt{\epsilon}S(B).
\end{align}
Thus, the corresponding inequality in (b) cannot hold for that $k$. Similarly, if the two inequalities in (b) hold for some $0\leq k\leq4$, then the corresponding inequality in (a) cannot be true.

Now, let us show that if $(a)$ holds for $k=0$, then it holds for all $0\leq k\leq4$. From the relations (\ref{4_elementi}) it follows
\begin{align}
b_{14}^{(k)}\pm b_{23}^{(k)} &= -\big{(}b_{14}^{(k-1)}\pm b_{23}^{(k-1)}\big{)}(\cos\phi_{k}\sin\psi_{k}\pm\sin\phi_{k}\cos\psi_{k}) \nonumber\\
&= -\big{(}b_{14}^{(k-1)}\pm b_{23}^{(k-1)}\big{)} \sin(\psi_{k}\pm\phi_{k}), \quad k\geq 1. \label{b14pmb23}
\end{align}
The relation (\ref{b14pmb23}) implies
\begin{equation}\label{slucaja}
|b_{14}^{(k)}\pm b_{23}^{(k)}| = |\sin(\psi_{k}\pm\phi_{k})|\cdots |\sin(\psi_{1}\pm\phi_{1})|\cdot|b_{14}^{(0)}\pm b_{23}^{(0)}| \leq |b_{14}^{(0)}\pm b_{23}^{(0)}|,
\end{equation}
for $k\geq1$. Therefore, if the two inequalities in (a) hold for $k=0$, then the relation (\ref{proofLemma3}) will hold for $k=1,2,3,4,5$. For any of these $k$ the relation (\ref{proofLemma3a}) also holds, proving that the inequalities in (b) cannot hold. Hence we can conclude that the both inequalities in~(a) hold for $0\leq k\leq4$.

Similarly, if the two inequalities in (b) hold for $k=0$, they hold for all $0\leq k\leq4$ and then the inequalities in (a) do not hold.
\end{proof}

We continue to prove (\ref{tmT1}) under the assumption $(a)$. The case (b) will be addressed later.

\subsection{The case $|b_{14}+b_{23}|\leq\sqrt{2}\delta_{1}S(B)$}\label{case_a}

Let us see what can be concluded for the rotation angles. Using $\phi$, $\psi$ for $\phi_{1}$, $\psi_{1}$, respectively, from the relations (\ref{4_elementi}) one easily obtains
\begin{align}
\big{(}b_{13}^{(1)}\big{)}^2+\big{(}b_{24}^{(1)}\big{)}^2
& = (b_{14}^2+b_{23}^2)(c_{\phi}^2 c_{\psi}^2+s_{\phi}^2s_{\psi}^2) - 4b_{14}b_{23}c_{\phi}c_{\psi}s_{\phi}s_{\psi} \nonumber \\
& = (b_{14}^2+b_{23}^2)(c_{\phi}c_{\psi}+s_{\phi}s_{\psi})^2 - 2(b_{14}^2+b_{23}^2)c_{\phi}c_{\psi}s_{\phi}s_{\psi} -  4b_{14}b_{23}c_{\phi}c_{\psi}s_{\phi}s_{\psi} \nonumber \\
& = (b_{14}^2+b_{23}^2)\cos^2(\phi-\psi) - 2(b_{14}+b_{23})^2c_{\phi}c_{\psi}s_{\phi}s_{\psi}. \label{prva}
\end{align}
Hence,
\begin{align*}
(b_{14}^2+b_{23}^2)\cos^2(\phi-\psi) & = 2(b_{14}+b_{23})^2c_{\phi}c_{\psi}s_{\phi}s_{\psi} + \big{(}b_{13}^{(1)}\big{)}^2+\big{(}b_{24}^{(1)}\big{)}^2 \leq \frac{1}{2}(b_{14}+b_{23})^2 + \delta_1^2S^2(B).
\end{align*}
We used the definition of $\delta_{1}$ from (\ref{4x4delta}). In the same way, one obtains
\begin{equation}\label{a14+a23}
\big((b_{14}^{(k-1)})^{2}+(b_{23}^{(k-1)})^{2}\big)\cos^{2}(\phi_{k}-\psi_{k}) \leq\frac{1}{2}(b_{14}^{(k-1)}+b_{23}^{(k-1)})^{2} + \delta_{k}^{2}S^{2}(B),
\end{equation}
for $k\geq1$. Using the relations (\ref{new1}), (\ref{4x4off}), the assumption (a) and Lemma~\ref{tm:lema3}, we conclude that the relation (\ref{a14+a23}) implies
$$(1-\epsilon)^2 S^2(B)\cos^2\big{(}\phi_{k}-\psi_{k}\big{)}  \leq S^2(B^{(k)})\cos^2\big{(}\phi_{k}-\psi_{k}\big{)} \leq (\delta_{k}^2+\delta_{k}^2)S^2(B),$$
for $1\leq k\leq5$. Thus
\begin{equation}\label{veliki_kutevi}
\cos\big(\phi_{k}-\psi_{k}\big) \leq\frac{\sqrt{2}}{1-\epsilon}\delta_{k} \leq1.4143\delta_{k}, \quad 1\leq k\leq5.
\end{equation}
Here we used $\epsilon=10^{-5}$. In (\ref{veliki_kutevi}) we have the strict inequalities when $\delta_{k}>0$ and that is certainly true for $1\leq k\leq4$.

\begin{Lemma}\label{tm:lemma4}
For the angles $\phi_{k}$, $\psi_{k}$, $1\leq k\leq5$, we have the following relations.
\begin{itemize}
\item[(i)] One of the following two relations holds
$$\phi_{k} = \frac{\pi}{4}-\alpha_k, \quad \psi_{k} =-\frac{\pi}{4}+\beta_k, \quad \alpha_k + \beta_k \leq 2.222\delta_{k}, \ \alpha_k\geq 0, \ \beta_k\geq 0,$$
$$\phi_{k} = -\frac{\pi}{4}+\alpha'_k, \quad \psi_{k} =\frac{\pi}{4}-\beta'_k, \quad \alpha'_k+\beta'_k \leq 2.222\delta_{k}, \ \alpha'_k\geq 0, \ \beta'_k\geq 0.$$
\item[(ii)] $|\sin(\phi_{k}+\psi_{k})| \leq |\phi_{k}+\psi_{k}| \leq 2.222\delta_{k}$.
\item[(iii)] $|\tan\phi_{k}+\tan\psi_{k}| \leq 4.444\delta_{k}$.
\item[(iv)] If $t\in\{|\tan\phi_{k}|,\tan\psi_k|\}$, then
$$0.98 < 1-4.444\delta_{k}\leq t\leq1 \quad \text{and} \quad 2\leq t+t^{-1} \leq 2+20.153\delta^2_{k}.$$
\item[(v)] $\max\left\{|\cot2\phi_{k}|,|\cot2\psi_{k}|,|\cot2\phi_{k}+\cot2\psi_{k}|\right\} \leq 4.49 \delta_{k}$.
\item[(vi)] $|\cot2\phi_{k}+\tan\phi_k+\cot2\psi_{k}+\tan\psi_k| \leq 10.08\delta_k^2 \leq 0.0454\delta_k$.
\item[(vii)] $2 \leq |\cot2\phi_{k}+\tan\phi_k-\cot2\psi_{k}-\tan\psi_k| \leq 2+20.15214\delta_k^2 \leq 2+0.0907\delta_k$.
\end{itemize}
\end{Lemma}

\begin{proof} The proof is technical and has been moved to~\ref{appB}.
\end{proof}

From the proof, one can easily check that in the assertions of Lemma~\ref{tm:lemma4}, the inequality signs $\leq$ and $\geq$ standing left to $\delta_{k}$ can be replaced by $<$ and $>$ respectively, provided that $\delta_{k}>0$ (which is true for $0\leq k\leq4$).

Let
\begin{equation}\label{nu}
\nu_k = \nu_k^{+}=\frac{|b_{14}^{(k)} + b_{23}^{(k)}|}{S(B)}, \qquad \nu_k^{-}=\frac{|b_{14}^{(k)} - b_{23}^{(k)}|}{S(B)}, \quad k\geq0; \quad \nu=\nu_0.
\end{equation}
Lemma~\ref{tm:lema3}(a) implies $\nu_{k}\leq\sqrt{2}\delta_{k+1}$, $0\leq k\leq4$.
From the relation (\ref{b14pmb23}) it follows
$$\nu_k=|\sin (\phi_k+\psi_k)|\nu_{k-1} \leq \nu_{k-1}, \quad \nu_k^{-}=|\sin (\phi_k-\psi_k)|\nu_{k-1}^{-} \leq \nu_{k-1}^{-}, \quad k\geq 0.$$
Hence by Lemma~\ref{tm:lemma4}(ii) we have
\begin{equation}\label{nu_k}
\nu_{k}\leq2.222\delta_{k}\nu_{k-1}, \quad 1\leq k\leq5.
\end{equation}

The next step is bounding $\delta_{k}$, $k=0,1,2$, by simple functions of the subsequent $\delta_{k}$s.

\begin{Lemma}\label{tm:lemma5}
The quantities $\delta_{k}$ satisfy the following inequalities
\begin{align*}
\delta_0 &< 9.348\delta_1\delta_2+19.371\delta_2\delta_3 + 9.646\delta_3\delta_4, \\
\delta_1 &< \delta_3 (9.464\delta_2 +7.874\delta_4) < 0.0388\delta_3, \\
\delta_2 &< \delta_4 (9.35\delta_3+7.778\delta_5) < 0.0771\delta_4.
\end{align*}
Hence we obtain
\begin{equation}\label{ocjened}
\delta_0 < 22.38\epsilon, \qquad \delta_1 < 17.5\epsilon, \qquad \delta_2 < 34.7\epsilon.
\end{equation}
\end{Lemma}

\begin{proof} The proof has been moved to~\ref{appB}.
\end{proof}

\begin{Lemma}\label{tm:lemma7}
For the pivot elements of $B^{(k)}$ we have:
\begin{itemize}
\item[(i)] $|b_{13}^{(k)}-b_{24}^{(k)}| \leq \nu_{k-1}S(B) \leq 2.222^{k-1}\delta_{k-1}\delta_{k-2}\cdots\delta_2\delta_1\sqrt{2}\delta_1S(B), \ 2\leq k\leq 6$. \\

In particular,
\begin{equation}\label{diffpivel}
\frac{|b_{13}^{(2)}-b_{24}^{(2)}|}{S(B)} < 3.15\delta_1^2, \quad
\frac{|b_{13}^{(3)}-b_{24}^{(3)}|}{S(B)} < 6.99\delta_2\delta_1^2, \quad
\frac{|b_{13}^{(4)}-b_{24}^{(4)}|}{S(B)} < 15.52\delta_3\delta_2\delta_1^2.
\end{equation}
\item[(ii)] $\big{(}b_{13}^{(k)}+b_{24}^{(k)}\big{)}^2 = 2\delta_k^2S^2(B)-\big{(}b_{13}^{(k)}-b_{24}^{(k)}\big{)}^2, \ k\geq0.$ \\

Hence,
$$|b_{13}^{(k)}+b_{24}^{(k)}| \leq \sqrt{2}\delta_kS(B), \quad k\geq 0,$$
and in particular
\begin{equation}\label{a13+a24-3}
\frac{|b_{13}^{(3)}+b_{24}^{(3)}|}{S(B)} > 1.41421\delta_3, \qquad \frac{|b_{13}^{(4)}+b_{24}^{(4)}|}{S(B)} > 1.41421\delta_4.
\end{equation}
\item[(iii)] $|\cot2\phi_{k}-\cot2\psi_{k}| \geq 1.999894\delta_k, \quad 3\leq k\leq4$.
\end{itemize}
\end{Lemma}

\begin{proof} The proof is technical and has been moved to~\ref{appB}.
\end{proof}

We come to the main part of the proof. So far, we derived the restrictions on the angles and other quantities, which are expressed in the previous lemmas. The question arises whether the diagonal elements, which enter into the definition of the angles, can allow all those limitations.

Consider the quantities
$$b^{(k)} = b_{11}^{(k)}-b_{22}^{(k)}-b_{33}^{(k)}+b_{44}^{(k)}, \quad k\geq 0.$$
Each $b^{(k)}$ will be expressed in two ways. On the one hand, we use (\ref{4_elementi}) and (\ref{4_ctg}) to obtain
\begin{align*}
b^{(k)}  & = b_{11}^{(k-1)}-b_{33}^{(k-1)}+2b_{13}^{(k-1)}\tan\phi_{k} + b_{22}^{(k-1)}-b_{44}^{(k-1)}+2b_{24}^{(k-1)}\tan\psi_{k} \\
& = 2b_{13}^{(k-1)}\cot2\phi_{k}+2b_{13}^{(k-1)}\tan\phi_{k}+2b_{24}^{(k-1)}\cot2\psi_{k}+
2b_{24}^{(k-1)}\tan\psi_{k} \\
& = 2b_{13}^{(k-1)}\big{(}\cot2\phi_{k}+\tan\phi_{k}\big{)}+2b_{24}^{(k-1)}\big{(}\cot2\psi_{k}+
\tan\psi_{k}\big{)} \\
& = \big{(}b_{13}^{(k-1)}+b_{24}^{(k-1)}\big{)}\big{(}\cot2\phi_{k}+\tan\phi_{k}+\cot2\psi_{k}+\tan\psi_{k}\big{)} + \\
& \quad + \big{(}b_{13}^{(k-1)}-b_{24}^{(k-1)}\big{)}\big{(}\cot2\phi_{k}+\tan\phi_{k}-\cot2\psi_{k}-\tan\psi_{k}\big{)}.
\end{align*}
Hence, by the assertions (vi) and (vii) of Lemma~\ref{tm:lemma4}, for $1\leq k\leq5$ we have
\begin{equation}\label{b3_1}
|b^{(k)}| \leq10.08\delta_{k}^{2}|b_{13}^{(k-1)}+b_{24}^{(k-1)}| + (2+20.15214\delta_{k}^{2})|b_{13}^{(k-1)}-b_{24}^{(k-1)}|.
\end{equation}
On the other hand, we use (\ref{4_ctg}) or (\ref{4_kutevi}) to obtain
\begin{align}
b^{(k)} & = (b_{11}^{(k)}-b_{33}^{(k)}) - (b_{22}^{(k)}-b_{44}^{(k)}) = 2b_{13}^{(k)}\cot2\phi_{k+1} - 2b_{24}^{(k)}\cot2\psi_{k+1} \nonumber\\
& = (b_{13}^{(k)} + b_{24}^{(k)})(\cot2\phi_{k+1}-\cot2\psi_{k+1}) + (b_{13}^{(k)}-b_{24}^{(k)})(\cot2\phi_{k+1}+\cot2\psi_{k+1}). \label{new4.47}
\end{align}
Therefore, we have
\begin{equation}\label{b3_2}
|(b_{13}^{(k)}+b_{24}^{(k)})(\cot2\phi_{k+1}-\cot2\psi_{k+1})| \leq |b_{13}^{(k)}-b_{24}^{(k)}|\cdot|\cot2\phi_{k+1}+\cot2\psi_{k+1}| + |b^{(k)}|.
\end{equation}
Using Lemma~\ref{tm:lemma7}(iii) we obtain
$$|(b_{13}^{(k)}+b_{24}^{(k)})(\cot2\phi_{k+1}-\cot2\psi_{k+1})| > |b_{13}^{(k)}+b_{24}^{(k)}| 1.999894\delta_{k+1}, \quad k=2,3.$$
Furthermore, using Lemma~\ref{tm:lemma4}(v), (\ref{b3_1}) and Lemma~\ref{tm:lemma7} we can bound from the above the right hand side of the inequality (\ref{b3_2}) divided by $S(B)$. We obtain
\begin{align*}
10.08\delta_k^2\sqrt{2}\delta_{k-1} + (2+20.15214\delta_k^2 + 4.49\delta_{k+1}2.222\delta_{k-1})2.222^{k-2}\delta_{k-2}\cdots\delta_1\sqrt{2}\delta_1 \\
\leq 14.256\delta_k^2\delta_{k-1} + (2+20.15214\delta_k^2+9.977\delta_{k+1}\delta_{k-1})2.222^{k-2}\delta_{k-2}\cdots\delta_1\sqrt{2}\delta_1.
\end{align*}
These inequalities hold for $2\leq k\leq4$. Hence, for $k=2,3$ we have
\begin{align}
&\frac{|b_{13}^{(k)}+b_{24}^{(k)}|}{S(B)}1.999894\delta_{k+1} \nonumber \\
& \leq (2+20.15214\delta_k^2+9.977\delta_{k+1}\delta_{k-1})
\cdot 2.222^{k-2}\delta_{k-2}\cdots\delta_1\sqrt{2}\delta_1 + 14.256\delta_k^2\delta_{k-1}. \label{k3case}
\end{align}
By inspecting the two cases, one checks that the case $k=3$ yields the contradiction with the relation~\eqref{4t} from the beginning of this proof. This is exactly what we need to prove the assertion~\eqref{tmT} of the theorem.
For $k=3$ the relation (\ref{k3case}) becomes
$$\frac{|b_{13}^{(3)}+b_{24}^{(3)}|}{S(B)}1.999894\delta_{4} \leq (2+20.15214\delta_3^2+9.977\delta_{4}\delta_{2})\cdot2.222\sqrt{2}\delta_{1}^2 + 14.256\delta_3^2\delta_{2}.$$
Using Lemma \ref{tm:lemma7}(ii) (actually the bound from (\ref{a13+a24-3})) the above relation implies
\begin{equation}\label{4.49}
2.828\delta_3\delta_4 < (6.285\delta_1+63.33\delta_1\delta_3^2+31.351\delta_4\delta_2\delta_1)\delta_1+14.256\delta_2\delta_3^2.
\end{equation}
From Lemma~\ref{tm:lemma5} we have $\delta_{1}< 0.0388\delta_{3}$ and $\delta_{2} < 0.0771\delta_{4}$. Moreover,
\begin{equation}\label{4.49a}
\delta_1 < \delta_3(9.464\delta_2 +7.874\delta_4) < \delta_3(9.464\cdot0.0771\delta_4 +7.874\delta_4) < 8.61\delta_3\delta_4.
\end{equation}
Dividing the inequality (\ref{4.49}) by $\delta_{3}\delta_{4}$ and using (\ref{4.49a}) and (\ref{ocjened}) we get the contradiction
\begin{align*}
2.828 & < (6.285\delta_1+63.33\delta_1\delta_3^2+31.351\delta_4\delta_2\delta_1)8.61+14.256\cdot0.0771\delta_3 \\
& < (6.285\cdot17.5\epsilon+63.33\cdot17.5\epsilon\cdot2\epsilon+31.351\cdot0.0045\cdot34.7\cdot17.5\epsilon^2)8.61 \\
& \qquad +14.256\cdot0.0771\cdot0.0045 < 0.014419.
\end{align*}
We used the bound $0.0045$ for $\delta_{3}$ and $\delta_{4}$, and $2\epsilon$ for $\delta_{3}^{2}$.

\subsection{The case $|b_{14}-b_{23}|\leq\sqrt{2}\delta_{1}S(B)$}\label{case_b}

The proof is similar as in the case $|b_{14}+b_{23}|\leq\sqrt{2}\delta_{1}S(B)$. We follow the lines of the proof above and modify it where necessary.

The expression on the right-hand side in the relation (\ref{prva}) can easily be brought to different form. We obtain
$$\big{(}b_{13}^{(1)}\big{)}^2+\big{(}b_{24}^{(1)}\big{)}^2 = (b_{14}^2+b_{23}^2)\cos^2(\phi+\psi) + 2(b_{14}-b_{23})^2c_{\phi}c_{\psi}s_{\phi}s_{\psi}.$$
Then the relation (\ref{a14+a23}) becomes
$$\left(\big{(}b_{14}^{(k-1)}\big{)}^2+\big{(}b_{23}^{(k-1)}\big{)}^2\right) \cos^2(\phi_{k}+\psi_{k}) \leq
\frac{1}{2}\big{(}b_{14}^{(k-1)}-b_{23}^{(k-1)}\big{)}^2 + \delta_{k}^2S^2(B), \quad k\geq1.$$
This implies
$$\cos(\phi_{k}+\psi_{k}) \leq \frac{\sqrt{2}}{1-\epsilon}\delta_k\leq1.4143\delta_k, \quad 1\leq k\leq5.$$
Lemma~\ref{tm:lemma4} has to be modified and we formulate it as a new lemma.

\begin{Lemma}\label{tm:lemma4_b}
For the angles $\phi_{k}$, $\psi_{k}$, $1\leq k\leq5$, we have the following relations.
\begin{itemize}
\item[(i)] One of the following two relations holds
$$\phi_{k}=\frac{\pi}{4}-\alpha_k, \quad \psi_{k}=\frac{\pi}{4}-\beta_k, \quad \alpha_k+\beta_k\leq2.222\delta_{k}, \ \alpha_k\geq0, \ \beta_k\geq0,$$
$$\phi_{k}=-\frac{\pi}{4}+\alpha'_k, \quad \psi_{k}=-\frac{\pi}{4}+\beta'_k, \quad \alpha'_k+\beta'_k\leq2.222\delta_{k}, \ \alpha'_k\geq 0, \ \beta'_k\geq0.$$
\item[(ii)] $|\sin(\phi_{k}-\psi_{k})| \leq |\phi_{k}-\psi_{k}|\leq2.222\delta_{k}$.
\item[(iii)] $|\tan\phi_{k}-\tan\psi_{k}| \leq 4.444\delta_{k}$.
\item[(iv)] If $t\in\{|\tan\phi_{k}|,\tan\psi_k|\}$, then
$$0.98 < 1-4.444\delta_{k}\leq t\leq1 \quad \text{and} \quad 2\leq t+t^{-1} \leq 2+20.153\delta^2_{k}.$$
\item[(v)] $\max \left\{|\cot2\phi_{k}|,|\cot2\psi_{k}|,|\cot2\phi_{k}-\cot2\psi_{k}|\right\} \leq 4.49 \delta_{k}$.
\item[(vi)] $|\cot2\phi_{k}+\tan\phi_k-\cot2\psi_{k}-\tan\psi_k|\leq 10.08\delta_k^2\leq0.0454\delta_k$.
\item[(vii)] $2 \leq |\cot2\phi_{k}+\tan\phi_k+\cot2\psi_{k}+\tan\psi_k| \leq2+20.15214\delta_k^2\leq2+0.0907\delta_k$.
\end{itemize}
\end{Lemma}

\begin{proof}
The proofs of these assertions are very similar or identical to the proofs of the corresponding assertions of Lemma~\ref{tm:lemma4}.
\end{proof}

Instead of $\nu_{k}$, we work with $\nu_{k}^{-}$. The relation (\ref{slucaja}) and the assertion (ii) of Lemma~\ref{tm:lemma4_b} imply
\begin{equation}\label{nu_k-}
\nu_{k}^{-} \leq 2.222\delta_{k}\nu_{k-1}^{-}, \quad 1\leq k\leq5.
\end{equation}
The statement of Lemma~\ref{tm:lemma5} does not have to be modified, but the proof needs minor changes. We have explained those changes in \ref{appB} under the title
``Proof of Lemma~\ref{tm:lemma5} in the case $|b_{14}-b_{23}|\leq\sqrt{2}\delta_{1}S(B)$''.

Lemma~\ref{tm:lemma7} has to be modified.

\begin{Lemma}\label{tm:lemma7_b}
For the pivot elements of $B^{(k)}$ we have:
\begin{itemize}
\item[(i)] $|b_{13}^{(k)}+b_{24}^{(k)}| \leq \nu_{k-1}^{-}S(B)\leq2.222^{k-1}\delta_{k-1}\delta_{k-2}\ldots\delta_2\delta_1\cdot \sqrt{2}\delta_1S(B), \ 2\leq k\leq6$. \\

In particular,
$$\frac{|b_{13}^{(2)}+b_{24}^{(2)}|}{S(B)} < 3.15\delta_1^2, \quad
\frac{|b_{13}^{(3)}+b_{24}^{(3)}|}{S(B)} < 6.99\delta_2\delta_1^2, \quad
\frac{|b_{13}^{(4)}+b_{24}^{(4)}|}{S(B)} < 15.52\delta_3\delta_2\delta_1^2.$$
\item[(ii)] $\big{(}b_{13}^{(k)}-b_{24}^{(k)}\big{)}^2 = 2\delta_k^2S^2(B)-\big{(}b_{13}^{(k)}+b_{24}^{(k)}\big{)}^2, \ \text{for} \ k\geq0$. \\

Hence,
\begin{equation}\label{4_zadnja}
|b_{13}^{(k)}-b_{24}^{(k)}| \leq \sqrt{2}\delta_k S(B), \quad k\geq 0,
\end{equation}
and in particular
$$\frac{|b_{13}^{(3)}-b_{24}^{(3)}|}{S(B)} > 1.41421\delta_3, \qquad \frac{|b_{13}^{(4)}-b_{24}^{(4)}|}{S(B)} > 1.41421\delta_4.$$
\item[(iii)] $|\cot2\phi_{k}+\cot2\psi_{k}| \geq 1.999894\delta_k, \quad 3\leq k\leq 4$.
\end{itemize}
\end{Lemma}

\begin{proof}
The proof is similar to the proof of Lemma~\ref{tm:lemma7}. We have moved it to~\ref{appB}.
\end{proof}

To prove the main assertion (\ref{tmT1}) we use the same $b^{(k)}$ as earlier. The assertions (vi) and (vii) of Lemma~\ref{tm:lemma4_b} yield
\begin{equation}\label{bk_zadnja}
|b^{(k)}| \leq10.08\,\delta_{k}^{2}|b_{13}^{(k-1)}-b_{24}^{(k-1)}| + (2+20.15214\delta_{k}^{2})|b_{13}^{(k-1)}+b_{24}^{(k-1)}|.
\end{equation}
Using (\ref{new4.47}) we obtain
$$\left|b_{13}^{(k)}-b_{24}^{(k)}\right|\left|\cot2\phi_{k+1}+\cot2\psi_{k+1}\right| \leq |b_{13}^{(k)}+b_{24}^{(k)}|\cdot|\cot2\phi_{k+1}-\cot2\psi_{k+1}|+|b^{(k)}|.$$
Using Lemma~\ref{tm:lemma7_b}(iii), the left-hand side can be bounded from below by
$1.999894 \cdot \delta_{k+1}\left|b_{13}^{(k)}-b_{24}^{(k)}\right|$ and for the case $k=3$ one can use Lemma~\ref{tm:lemma7_b}(ii) to further reduce it to
$2.828 \delta_{3}\delta_{4} S(B)$.

Using (\ref{bk_zadnja}), (\ref{4_zadnja}), Lemma~\ref{tm:lemma4_b}(v), and Lemma~\ref{tm:lemma7_b}(i), the right-hand side divided by $S(B)$ can be bounded from above by
\begin{align*}
10.08\delta_k^2\sqrt{2}\delta_{k-1} + (2+20.15214\delta_k^2+4.49\delta_{k+1}2.222\delta_{k-1})2.222^{k-2}\delta_{k-2}\cdots\delta_1\sqrt{2}\delta_1 \\
\leq 14.256\delta_k^2\delta_{k-1} + (2+20.15214\delta_k^2+9.977\delta_{k+1}\delta_{k-1})2.222^{k-2}\delta_{k-2}\cdots\delta_1\sqrt{2}\delta_1.
\end{align*}
For $k=3$, after dividing by $S(B)$, one obtains
$$2.828 \delta_3\delta_4 <  (6.285\delta_1+63.33\delta_1\delta_3^2+31.351\delta_4\delta_2\delta_1)\delta_1+14.256\delta_2\delta_3^2.$$
which is the same inequality as (\ref{4.49}). The rest of the proof is the same as earlier. 
\qed

\bigskip

At this point we would like to make a few comments.
\begin{itemize}
\item Theorem~\ref{tm:main} obviously holds with somewhat larger $\epsilon$, e.g. one can try to complete the proof with $\epsilon=10^{-4}$. On the other hand, a small $\epsilon$ from the proof exposes the possibility of the very small reduction of $S(A)$ within one cycle. This happens when an underlaying matrix has a special structure. The next section deals with this issue.
\item Although $1-10^{-5}$ is a small decrease of the off-norm within two cycles, the result does not mean that the convergence of the method should be slow. The proof is concentrated on the worst case scenario. Typically, the slower the method is within one cycle, the faster it is in the next cycle. Example~\ref{ex:matrix1} indicates that behavior.
\item In this convergence proof we have explicitly used the diagonal elements of $A$, which is unusual when the reduction of $S(A)$ is considered. Usually, only the off-diagonal elements and the bounds on rotation angles are used (e.g. \cite{hen+zim-68,SS89,wil-62,har-91}). In that case the proof is valid for a more general iterative process used in the global convergence analysis of Jacobi-type processes which use nonorthogonal transformation matrices \cite{H15,B}.
\end{itemize}

\section{The slow off-norm reduction within one cycle}\label{sec5}

As it can be seen from the above theory, the decrease of the off-norm after one cycle of the Jacobi method under the strategy $I_{1}$ can be small. Here we give an example from \cite{B}, where the relative decrease of the off-norm after one cycle is less then $10^{-50}$.

\begin{Example}\label{ex:matrix1}
Let
$$H=H^{(0)}
=\left[\begin{array}{cccc}
1+p_{1}+p_{2} & 0 & \epsilon+p_{1} & -1+p_{1} \\[-1pt]
0 & 1+p_{2} & 1 & -\epsilon \\[-1pt]
\epsilon+p_{1} & 1 & 1+p_{1} & 0 \\[-1pt]
-1+p_{1} & \epsilon& 0 & 1 \\
\end{array}\right],
$$
with $\epsilon=10^{-52}$, $p_{1}=\epsilon$, $p_{2}=\epsilon\sqrt{\epsilon}$.

We have used MATLAB Symbolic Math Toolbox, in particular the Variable-precision arithmetic with $2l$ digits, to compute the matrix iterates under the cyclic Jacobi method defined by the strategy $I_{1}$. We display the off-norm of each iterate to $l$ significant digits. For $l=50$ we obtain
$$S(H)=1.41421356237309504880168872420969807856967187537694.$$
As we can see from the table below, during the first cycle the off-norm of $H$ does not change in the first $50$ decimal places. But later it drops rapidly, especially in the $8$th step.
\bigskip

\begin{small}
\begin{center}
\begin{tabular}{|c|c|c|}
  \hline
  $k$ & $(i(k),j(k))$ &  $S(H^{(k)})$ \\
  \hline
  $1$ & $(1,3)$ & $1.41421356237309504880168872420969807856967187537694$\hspace{5ex} \\
  $2$ & $(2,4)$ & $1.41421356237309504880168872420969807856967187537694$\hspace{5ex} \\
  $3$ & $(1,4)$ & $1.41421356237309504880168872420969807856967187537694$\hspace{5ex} \\
  $4$ & $(2,3)$ & $1.41421356237309504880168872420969807856967187537694$\hspace{5ex} \\
  $5$ & $(1,2)$ & $1.41421356237309504880168872420969807856967187537694$\hspace{5ex} \\
  $6$ & $(3,4)$ & $1.41421356237309504880168872420969807856967187537694$\hspace{5ex} \\ \hline
  $7$ & $(1,3)$ & $0.99999999999999999999999999999999999999999999999999$\hspace{5ex} \\
  $8$ & $(2,4)$ & $0.17677669529663688110021108266947024663734760219051{e{-}26}$ \\
  \hline
\end{tabular}
\end{center}
\end{small}
\end{Example}
\bigskip

In general, one can always find a matrix $A(\epsilon)$ such that the decrease of the off-norm after one cycle of the Jacobi method under the strategy $I_{1}$ is arbitrary small and depends only on $\epsilon$.

\begin{Proposition}\label{tm:prop5.1}
Let $0<\epsilon\leq10^{-5}$,
$$H(\epsilon)=\left[
    \begin{array}{cccc}
      \epsilon+\epsilon^{1.5} & 0 & 2\epsilon & -1+\epsilon \\
      0 & \epsilon^{1.5} & 1 & -\epsilon \\
      2\epsilon & 1 & \epsilon & 0 \\
      -1+\epsilon & -\epsilon & 0 & 0 \\
    \end{array}
  \right],$$
and let the cyclic Jacobi method defined by the strategy $I_{1}$ be applied to $H(\epsilon)$, thus generating the matrices $H^{(0)}=H(\epsilon),H^{(1)},\ldots{}$. After completing one full sweep we have
$$S^{2}(H^{(6)})> (1-17\epsilon)S^{2}(H^{(0)}).$$
\end{Proposition}

\begin{proof} The proof is lengthy and technical, so it has been moved to~\ref{appC}.
\end{proof}

We end the paper with the following important theorem.

\begin{Theorem}\label{tm:teo5.2}
For every $0<\epsilon<1$ and $n\geq4$, there exists a symmetric matrix $A(\epsilon)$ of order $n$, depending on $\epsilon$ and a cyclic strategy $I$, such that
\begin{equation}\label{pr_main}
S(A^{(N)})> (1-\epsilon)S(A^{(0)}), \quad N=\frac{n(n-1)}{2}.
\end{equation}
Here $A^{(0)}=A(\epsilon)$ and $A^{(N)}$ is obtained from $A^{(0)}$ by applying a full cycle of the Jacobi method under the strategy $I$.
\end{Theorem}

\begin{proof}
Let $n=4$, $0<\epsilon<10^{-5}$, $\epsilon'= (2\epsilon-\epsilon^{2})/17$ and $A(\epsilon) = H(\epsilon')$, where $H(\epsilon')$ is from Proposition~\ref{tm:prop5.1}.
Let $I=I_{1}=I_{\mathcal{O}_{1}}$.
Proposition~\ref{tm:prop5.1} yields to $A^{(0)}= A(\epsilon)$,
$$S^{2}(A^{(6)}) > (1-17\epsilon')S^{2}(A^{(0)}).$$
Since
$$\sqrt{1-17\epsilon'} = \sqrt{1-2\epsilon+\epsilon^{2}}=1-\epsilon,$$
the proof is completed in this case.

If $10^{-5}\leq\epsilon<1$, then $1-10^{-5} \geq1-\epsilon>0$.
Hence, we can choose $A(\epsilon)= H(\frac{2\cdot10^{-5}-10^{-10}}{17})$ to obtain $S(A^{(6)}) > (1-10^{-5})S(A) \geq(1-\epsilon)S(A)$.

Let $n>4$ and let
$$A^{(0)} = \left[\begin{array}{cc}
A_{11}^{[0]} & 0 \\
0 & A_{22}^{[0]}
\end{array}\right]$$
be a symmetric matrix of order $n$ with the following properties.
\begin{itemize}
\item[(i)] $A_{11}^{[0]}$ is of order $4$ such that $S(A_{11}^{(6)}) > (1-\epsilon)S(A_{11}^{[0]})$ holds when one full cycle of the Jacobi method under the strategy $I_{1}$ is applied to $A_{11}^{[0]}$. This follows from (\ref{pr_main}) because we have proved the theorem for $n=4$.
\item[(ii)] The block $A_{22}^{[0]}$ is diagonal.
\end{itemize}
The pivot strategy $I$ is defined by $I=I_{\mathcal{O}}$, $\mathcal{O} = [\mathcal{O}_{1},\mathcal{O}_{22},\mathcal{O}_{12}]$, where $\mathcal{O}_{22}$ is any ordering of the set
$\mathcal{S}_{22}=\{(5,6),\ldots,(5,n),\ldots,(n-2,n-1),(n-2,n),(n-1,n)\}$ and $\mathcal{O}_{12}$ is any ordering of the set $\mathcal{S}_{12}=\{(1,5),\ldots,(1,n),\ldots,(4,5),\ldots,(4,n)\}$.

Obviously, the whole sweep on $A(\epsilon)$ reduces to the sweep on $A_{11}^{[0]}$ under the strategy~$I_{\mathcal{O}_{1}}$ since all other Jacobi angles are zero.
\end{proof}

Let us show that the blocks $A_{12}^{[0]}$ and $A_{22}^{[0]}$ of the matrix $A^{(0)}= A(\epsilon)$ can be chosen such that all their entries are nonzero. Indeed, we can make other sets of the assumptions on~$A(\epsilon)$. One such set of the assumptions is the following.
\begin{itemize}
\item[(i)] $A_{11}^{[0]}$ is of order $4$ and such that
\begin{equation}\label{pr_main1a}
S(A_{11}^{(6)})> (1-\frac{\epsilon}{2})S(A_{11}^{[0]})
\end{equation}
holds when a full cycle of the Jacobi method under the strategy $I_{1}$ is applied to~$A_{11}^{[0]}$. The existence of such an $A_{11}^{[0]}$ follows from (\ref{pr_main}) because we have proved the theorem for $n=4$.
\item[(ii)]  We have $\|A_{12}^{[0]}\|_{F}\leq\epsilon^{\nu}S(A_{11}^{[0]})$, where $\nu$ satisfies $n\epsilon^{\nu-2} <1$.
\item[(iii)] We have $S(A_{22}^{[0]})\leq\epsilon^{2}S(A_{11}^{[0]})$ and
$$\diag (A_{22}^{[0]}) = \diag(\delta_1,\delta_2,\ldots,\delta_{n-4}),$$
where
\begin{equation}\label{pr_A22}
\delta_i \in \big{\{}t \ | \ t< -\|A_{11}^{[0]}\|_2-1\} \cup \{t \ | \ t > \|A_{11}^{[0]}\|_2+1\big{\}}.
\end{equation}
\end{itemize}
The pivot strategy $I$ is defined as in the proof of Theorem~\ref{tm:teo5.2}.

To keep the paper shorter we do not give a rigorous proof, but we make few essential remarks. After completing the sweep on $A_{22}^{[6]}$ the inequality (\ref{pr_main1a}) still holds, only the superscript $(6)$ on the left-hand side has to be replaced by $(M)$, $M=6+(n-4)(n-\nobreak 5)/2$. We also have $S(A_{22}^{(M)})\leq\epsilon^{2}S(A_{11}^{[0]})$ and $\|A_{12}^{(M)}\|_{F}=\|A_{12}^{[0]}\|_{F}\leq\epsilon^{\nu}S(A_{11}^{[0]})$. Due to the condition (\ref{pr_A22}) all later angles will be bounded by some
multiples of $\epsilon^{\nu}$. The sum of squares of the last $N-M$ pivot elements will be bounded by some multiple of
$n\epsilon^{2\nu}S^{2}(A_{11}^{[0]})< \epsilon^{4}S^{2}(A_{11}^{[0]})$, which will eventually yield the required result.

\section*{Acknowledgements}
The authors are thankful to the anonymous referees for their excellent remarks which improved the readability of the paper.

\appendix
\setcounter{equation}{0}

\section{Proofs related to Section~\ref{sec3}}\label{appA}

\subsection{Proof of Proposition~\ref{tm:prva}}

We prove the proposition for $P=P_{12}$. The proof for the case $P=P_{34}$ is similar. Since the both Jacobi processes are cyclic, it is sufficient to prove the proposition for $r=1,2,3$. Let
$$A^{(6)} =  U^TAU \qquad \mathsf{A}^{(6)} =  \mathsf{U}^T\mathsf{A}^{(0)}\mathsf{U} = \mathsf{U}^TP^TAP\mathsf{U},$$
where
\begin{align*}
U &= R(1,4,\theta^{(0)})R(2,3,\theta^{(1)})R(1,3,\theta^{(2)})R(2,4,\theta^{(3)})R(1,2,\theta^{(4)})R(3,4,\theta^{(5)}),\\
\mathsf{U} &= R(1,3,\phi^{(0)})R(2,4,\phi^{(1)})R(1,4,\phi^{(2)})R(2,3,\phi^{(3)})R(1,2,\phi^{(4)})(3,4,\phi^{(5)}).
\end{align*}
Note that $P=P^{T}$. Let us inspect the product $P\mathsf{U}=P^{T}\mathsf{U}$. We have
\begin{align*}
P\mathsf{U} &= P^TR(1,3,\phi^{(0)})P\cdot P^TR(2,4,\phi^{(1)})P\cdot P^TR(1,4,\phi^{(2)})P\cdot P^TR(2,3,\phi^{(3)})P \\
& \qquad \cdot P^TR(1,2,\phi^{(4)})P\cdot P^TR(3,4,\phi^{(5)})P\cdot P^T \\
&= R(2,3,\phi^{(0)})R(1,4,\phi^{(1)})R(2,4,\phi^{(2)})R(1,3,\phi^{(3)})R(1,2,-\phi^{(4)})R(3,4,\phi^{(5)})P^T \\
&= R(1,4,\phi^{(1)})R(2,3,\phi^{(0)})R(1,3,\phi^{(3)})R(2,4,\phi^{(2)})R(1,2,-\phi^{(4)})R(3,4,\phi^{(5)})P,
\end{align*}
It remains to show that $\theta^{(0)}=\phi^{(1)}$, $\theta^{(1)}=\phi^{(0)}$, $\theta^{(2)}=\phi^{(3)}$, $\theta^{(3)}=\phi^{(2)}$, $\theta^{(4)}=-\phi^{(4)}$, $\theta^{(5)}=\phi^{(5)}$. Since
\begin{equation}\label{obaA}
A^{(0)}=\left[
    \begin{array}{cccc}
      a_{11} & a_{12} & a_{13} & a_{14} \\
      a_{12} & a_{22} & a_{23} & a_{24} \\
      a_{13} & a_{23} & a_{33} & a_{34} \\
      a_{14} & a_{24} & a_{34} & a_{44} \\
    \end{array}
  \right], \qquad
 \mathsf{A}^{(0)}=\left[
    \begin{array}{cccc}
      a_{22} & a_{12} & a_{23} & a_{24} \\
      a_{12} & a_{11} & a_{13} & a_{14} \\
      a_{23} & a_{13} & a_{33} & a_{34} \\
      a_{24} & a_{14} & a_{34} & a_{44} \\
    \end{array}
  \right],
\end{equation}
it immediately follows from (\ref{kut1}) that $\theta^{(0)}=\phi^{(1)}$, $\theta^{(1)}=\phi^{(0)}$.
Thus, after completing the first two steps in each of the two processes, we have
$$\mathsf{A}^{(2)} =  P^TR_{23}^T (\theta^{(1)}) R_{14}^T(\theta^{(0)})A R_{14}(\theta^{(0)})R_{23}(\theta^{(1)})P = P^TA^{(2)}P.$$
This shows that the relation (\ref{obaA}) holds if $A^{(0)}$ and $\mathsf{A}^{(0)}$ are replaced by $A^{(2)}$ and $\mathsf{A}^{(2)}=P^{T}A^{(2)}P$, respectively.
Checking the angle formula (\ref{kut1}) we find that $\theta^{(2)}=\phi^{(3)}$, $\theta^{(3)}=\phi^{(2)}$ and therefore $\mathsf{A}^{(4)}=P^{T}A^{(4)}P$.
The last check is the easiest one since the denominators in (\ref{kut1}) for the angles $\theta^{(4)}$ and $\theta^{(5)}$ are opposite to those for the angles $\phi^{(4)}$ and $\phi^{(5)}$.
\qed

\subsection{Proof of Proposition~\ref{tm:prop1}}

Let us denote $B^{(k)}=\mathcal{T}^{k}(A)$, $k\geq0$. An easy calculation shows that
\begin{align}
Q & =[e_1 \ e_3 \ e_4 \ \mbox{-}e_2], \qquad Q^2 =[e_1 \ e_4 \ \mbox{-}e_2 \ \mbox{-}e_3], \qquad Q^3 =[e_1 \ \mbox{-}e_2 \ \mbox{-}e_3 \ \mbox{-}e_4], \label{Q_1}\\
Q^4 & = [e_1 \ \mbox{-}e_3 \ \mbox{-}e_4 \ e_2], \qquad Q^5  =[e_1 \ \mbox{-}e_4 \ e_2 \ e_3], \qquad Q^6 =[e_1 \ e_2 \ e_3 \ e_4].\label{Q_2}
\end{align}
Hence $B^{(6)}=\mathcal{T}^{6}(A)=(Q^{T})^{6}A^{(12)}Q^{6}=A^{(12)}$ and it is sufficient to show that the relation (\ref{Tprop}) holds for $0\leq k\leq6$. We shall show
$$B^{(k)}=(Q^{k})^T A^{(2k)}Q^{k}, \quad 0\leq k\leq 6.$$
Consider two processes, the first one is defined by the relation $B^{(k)}=\mathcal{T}^{k}(A)$, $k\geq0$, and the second one is the Jacobi method under the strategy $I_{1}$.
These two processes generate the matrices $B^{(k)}=(b_{rs}^{(k)})$, $k\geq0$, and $A^{(l)}=(a_{rs}^{(l)})$, $l\geq0$, respectively.
The rotation angles at the step $k$ of the first process will be denoted by $\phi_{k}$ and $\psi_{k}$, $k\geq1$. The rotation angle at the step $l$ of the Jacobi method will be denoted by $\varphi^{(l)}$, $l\geq0$.
Thus, $\phi_{k}$ and $\psi_{k}$ are used to compute $B^{(k)}$, while $\varphi^{(l-1)}$ is used to compute $A^{(l)}$.

For $k=0$ the assertion (\ref{Tprop}) takes the form $B^{(0)}=A^{(0)}$ which is correct since $B^{(0)}=A$ and $A^{(0)}=A$.

Let $k=1$. Then
$$B^{(1)} = \mathcal{T}(A) = Q^TR(2,4,\psi_1)^TR(1,3,\phi_1)^TAR(1,3,\phi_1)R(2,4,\psi_1)Q.$$
By Definition~\ref{T} angles $\phi_{1}$ and $\psi_{1}$ are the Jacobi angles which annihilate the elements of~$A$ at positions $(1,3)$ and $(2,4)$. Therefore, we have $\phi_{1}=\varphi^{(0)}$ and $\psi_{1}=\varphi^{(1)}$, and consequently $A^{(2)}=R(2,4,\psi_{1})^{T}R(1,3,\phi_{1})^{T}AR(1,3,\phi_{1})R(2,4,\psi_{1})$.
Thus, $B^{(1)} = Q^{T}A^{(2)}Q$, which had to be proved.

Let $k=2$. We use the fact that the assertion (\ref{Tprop}) holds for $k=1$. Using the relations (\ref{Q_1}) and (\ref{Q'XQ}) one obtains
\begin{align*}
B^{(2)} & =\mathcal{T}(B^{(1)}) = \mathcal{T}(Q^TA^{(2)}Q) = (R(1,3,\phi_2)R(2,4,\psi_2)Q)^TQ^TA^{(2)}QR(1,3,\phi_2)R(2,4,\psi_2)Q \\
& = Q^TR(2,4,\psi_2)^TR(1,3,\phi_2)^T{\scriptsize \left[
    \begin{array}{cccc}
      a_{11}^{(2)} & 0 & a_{14}^{(2)} & -a_{12}^{(2)} \\
      0 & a_{33}^{(2)} & a_{34}^{(2)} & -a_{23}^{(2)} \\
      a_{14}^{(2)} & a_{34}^{(2)} & a_{44}^{(2)} & 0 \\
      -a_{12}^{(2)} & -a_{23}^{(2)} & 0 & a_{22}^{(2)} \\
    \end{array}
  \right]}R(1,3,\phi_2)R(2,4,\psi_2)Q.
\end{align*}
For the rotation angles which annihilate the elements at positions $(1,3)$ and $(2,4)$ we have
$$\tan2\phi_2=\frac{2a_{14}^{(2)}}{a_{11}^{(2)}-a_{44}^{(2)}}, \qquad
\tan2\psi_2=\frac{-2a_{23}^{(2)}}{a_{33}^{(2)}-a_{22}^{(2)}} = \frac{2a_{23}^{(2)}}{a_{22}^{(2)}-a_{33}^{(2)}},$$
hence,
\begin{equation}\label{Tkutevi2}
\varphi^{(2)}=\phi_{2},\qquad\varphi^{(3)}=\psi_{2}.
\end{equation}
The relation $B^{(2)}=(Q^{2})^{T} A^{(4)}Q^{2}$ will hold provided that
\begin{equation}\label{QR2}
QR(1,3,\phi_{2})R(2,4,\psi_{2})Q = R(1,4,\varphi^{(2)})R(2,3,\varphi^{(3)})Q^{2},
\end{equation}
and the relation (\ref{QR2}) will hold provided that
\begin{align}
R(1,4,\varphi^{(2)}) & = QR(1,3,\phi_2)Q^T, \label{Q2_1} \\
R(2,3,\varphi^{(3)}) & = QR(2,4,\psi_2)Q^T. \label{Q2_2}
\end{align}
It is easy to see that the relations (\ref{Q2_1}) and (\ref{Q2_2}) follow from the relations (\ref{Q'XQ}) and (\ref{Tkutevi2}).

The proof for $k=3,4,5,6$ proceeds in the same manner as for $k=2$, but with different indices.
\qed

\subsection{Proof of Lemma~\ref{4x4lema}}

We shall use the notation from the proof of Proposition~\ref{tm:prop1}.
\begin{itemize}
\item[(i)] If $a_{14}=0$ and $a_{23}=0$, then $S^{2}(A)=a_{13}^{2}+a_{24}^{2}$ and
$$S^{2}(A^{(2)}) = S^{2}(A)-(a_{13}^{2}+a_{24}^{2}) = 0.$$
\item[(ii)] Since the first two pivot elements $a_{13}$ and $a_{24}$ are zero, the corresponding rotation angles $\varphi^{(0)}$ and $\varphi^{(1)}$ are zero as well, and $A^{(2)}=A$. The next two pivot elements $a_{14}$ and $a_{23}$ are the only possibly nonzero off-diagonal elements. Hence,
    $S^{2}(A^{(4)})=S^{2}(A)-a_{14}^{2}+a_{23}^{2}=0$.
\item[(iii)] Since $a_{13}=0$, we have $\phi_{1}=0$. The relations (\ref{4_elementi}) imply
$$(b_{13}^{(1)})^2+(b_{24}^{(1)})^2 = (a_{14}^2+a_{23}^2)\cos^2\psi_1 \geq \frac{1}{2}(a_{14}^2+a_{23}^2).$$
We used the assumption $\psi_{1}\in[-\frac{\pi}{4},\frac{\pi}{4}]$. Using the relation (\ref{STA}), we have
\begin{align*}
S^2(A^{(4)}) & = S^2(B^{(2)}) = S^2(B^{(1)})-\big{(}(b_{13}^{(1)})^2+(b_{24}^{(1)})^2\big{)} \leq \big{(}S^2(A)-a_{24}^2\big{)} - \frac{1}{2}(a_{14}^2+a_{23}^2) \\
& \leq S^2(A) - \frac{1}{2}(a_{24}^2+a_{14}^2+a_{23}^2) = S^2(A)-\frac{1}{2}S^2(A) = \frac{1}{2}S^2(A).
\end{align*}
\item[(iv)] The proof is same as (iii), only $a_{13}$ and $\psi_{1}$ are used instead of $a_{24}$ and $\phi_{1}$.
\item[(v)] Let $a_{14}=0$. From the relations (\ref{4_elementi}) we get
$$(b_{13}^{(1)})^2+(b_{24}^{(1)})^2 = a_{23}^2 (\cos^2\phi_1\cos^2\psi_1+\sin^2\phi_1\sin^2\psi_1) \geq \frac{1}{2}a_{23}^2.$$
We bounded the expression
$\cos^{2}\phi_{1}\cos^{2}\psi_{1}+\sin^{2}\phi_{1}\sin^{2}\psi_{1}$
using the function $f(x,y)= 1-(x^{2}+y^{2})+2x^{2}y^{2}$ for $x=\sin\phi_{1}$, $y=\sin\psi_{1}$
on $[-\frac{\sqrt{2}}{2}, \frac{\sqrt{2}}{2}]\times[-\frac{\sqrt{2}}{2}, \frac{\sqrt{2}}{2}]$.
The minimum of that function equals $\frac{1}{2}$. Hence,
\begin{align*}
S^2(A^{(4)}) & = S^2(B^{(2)}) = S^2(B^{(1)})-\big{(}(b_{13}^{(1)})^2+(b_{24}^{(1)})^2\big{)} \leq \big{(}S^2(A)-a_{13}^2-a_{24}^2\big{)} - \frac{1}{2}a_{23}^2 \\
& \leq S^2(A) - \frac{1}{2}(a_{13}^2+a_{24}^2+a_{23}^2) = S^2(A)-\frac{1}{2}S^2(A) = \frac{1}{2}S^2(A).
\end{align*}
\item[(vi)] The proof is same as $(v)$, only $a_{14}$ is used instead of $a_{23}$.
\qed
\end{itemize}

\section{Proofs related to Section~\ref{proof}}\label{appB}

\subsection{Proof of Lemma~\ref{tm:lemma4}}

First, note that the following two inequalities hold:
\begin{align}
\label{tm:lema2i} |x| & \geq\frac{\pi}{2}(1-\cos(x)), \quad x\in\left[ -\frac{\pi}{2},\frac{\pi}{2}\right], \\
\label{tm:lema2ii} \left| \tan(x_{1})-\tan(x_{2}) \right| & \leq2|x_{1} -x_{2}|, \quad x_{1},x_{2}\in\left[-\frac{\pi}{4},\frac{\pi}{4}\right].
\end{align}

\begin{itemize}
\item[(i)] Relations \eqref{veliki_kutevi} and \eqref{tm:lema2i} imply
$$|\phi_{k}-\psi_{k}| \geq (1-1.4143\delta_{k})\frac{\pi}{2} \geq \frac{\pi}{2}-2.222\delta_{k}, \quad 1\leq k\leq 5.$$
The assertion follows from the fact that the rotation angles are from the interval $\displaystyle[-\frac{\pi}{4},\frac{\pi}{4}]$.
\item[(ii)]  The assertion follows from (i).
\item[(iii)] Using the relations \eqref{tm:lema2ii} and \eqref{tm:lema2i} we have
$$|\tan\phi_{k}+\tan\psi_{k}| = |\tan|\phi_{k}|-\tan |\psi_{k}|| \leq 2\max\left\{|\alpha_k-\beta_k|,|\alpha_k'-\beta_k'|\right\} \leq 4.444 \delta_{k}.$$
\item[(iv)]  From the relations \eqref{tm:lema2ii} and \eqref{tm:lema2i} we have either
$$|\tan\phi_{k}-1| \leq 2\left|\phi_{k}-\frac{\pi}{4}\right|\leq 4.444\delta_{k}, \qquad |\tan\psi_{k}+1|\leq2\left|\psi_{k}+\frac{\pi}{4}\right|\leq 4.444\delta_{k},$$
or
$$|\tan\phi_{k}+1| \leq2\left|\phi_{k}+\frac{\pi}{4}\right| \leq 4.444\delta_{k}, \qquad |\tan\psi_{k}-1|\leq2\left|\psi_{k}-\frac{\pi}{4}\right| \leq 4.444\delta_{k}.$$
Hence, using (\ref{new3}) one obtains the lower bound for the tangents. The upper bound is obvious since the angles lie in the segment $[-\pi/4,\pi/4]$.
For the latter assertion note that $x+x^{-1}\geq2$ holds for any real $x\ne0$ and equality is attained only for $x=1$. Recall that we can write $t=1-\gamma$, $0\leq\gamma\leq4.444\delta_{k}$. Hence,
\begin{align*}
t+\frac{1}{t} & = 1-\gamma + \frac{1}{1-\gamma}= 1-\gamma + 1+\gamma +\frac{\gamma^2}{1-\gamma} \\
& \leq 2+\frac{4.444^2}{1-4.444\cdot 0.0045}\delta_k^2 \leq 2+20.15214\delta_k^2.
\end{align*}
\item[(v)]  Specifying $t=\min\{|\tan\phi_{k}|,|\tan\psi_{k}|\}$ we have
\begin{align*}
\max\left\{|\cot2\phi_{k}|,|\cot2\psi_{k}|\right\} & = \frac{1-t^2}{2t} =
\frac{1}{2}\left|\frac{1}{t} -t \right| = \gamma + \frac{1}{2}\frac{\gamma^2}{1-\gamma} \\
& \leq \left(1+\frac{1}{2}\frac{4.444\cdot 0.0045}{1-4.444\cdot 0.0045}\right)4.444\delta_k \leq 4.48935\delta_k.
\end{align*}
Since $\cot2\phi_{k}$ and $\cot2\psi_{k}$ have the opposite sign, the absolute value of their sum cannot be larger than the larger term.
\item[(vi)] Let $\eta_{k} = \cot2\phi_{k} + \tan\phi_{k} + \cot2\psi_{k} + \tan\psi_{k}$.
Using the notation and ideas from the proof of (iv), we have
$$\eta_k = \frac{1}{2}\big{(} t_{\phi_k}+\frac{1}{t_{\phi_k}} +t_{\psi_k}+\frac{1}{t_{\psi_k}}\big{)} =
\frac{1}{2}\sigma_k\big{(}\frac{\gamma_{\phi_k}^2}{1-\gamma_{\phi_k}} - \frac{\gamma_{\psi_k}^2}{1-\gamma_{\psi_k}} \big{)}.$$
This implies
$$|\eta_k | \leq \frac{1}{2}\max\big{\{}\frac{\gamma_{\phi_k}^2}{1-\gamma_{\phi_k}}\frac{\gamma_{\psi_k}^2}{1-\gamma_{\psi_k}}\big{\}}
\leq 10.07607\delta_k^2 \leq 0.045343\delta_k.$$
\item[(vii)] The proof is similar to the proof of (vii). If
$$|\mu_k|=\frac{1}{2}\left|2+\frac{\gamma_{\phi_k}^2}{1-\gamma_{\phi_k}} +2+\frac{\gamma_{\psi_k}^2}{1-\gamma_{\psi_k}}\right| \leq 2+20.15214\delta_k^2 \leq 2+0.0907\delta_k.$$
\qed
\end{itemize}

\subsection{Proof of Lemma~\ref{tm:lemma5}}

In terms of the elements of matrices $B^{(k)}$ and $B^{(k-1)}$ for $k\geq1$ the angle formulas (\ref{4_kutevi}) take the form
\begin{equation}\label{4_ctg}
\frac{b_{11}^{(k-1)}-b_{33}^{(k-1)}}{2b_{13}^{(k-1)}}=\cot\big{(}2\phi_{k}\big{)}, \qquad \frac{b_{22}^{(k-1)}-b_{44}^{(k-1)}}{2b_{24}^{(k-1)}}=\cot\big{(}2\psi_k\big{)}.
\end{equation}
The relation (\ref{nu_k}) and Lemma~\ref{tm:lema3}(a) imply
\begin{equation}\label{veliki_elementi_plus}
\nu_k \leq 2.222^k\delta_k\cdots\delta_1\nu \leq 2.222^k\delta_k\cdots\delta_2\delta_1\cdot\sqrt{2}\delta_1 S(B), \quad 1\leq k\leq 5.
\end{equation}

From the relations (\ref{4_elementi}) for $k\geq1$ we have
\begin{align*}
b_{11}^{(k)}-b_{22}^{(k)} & = b_{11}^{(k-1)}-b_{33}^{(k-1)}+2b_{13}^{(k-1)}\tan\phi_{k}, \\
b_{44}^{(k)}-b_{33}^{(k)} & = b_{22}^{(k-1)}-b_{44}^{(k-1)}+2b_{24}^{(k-1)}\tan\psi_{k}.
\end{align*}
Combining that with the angle formulas (\ref{4_ctg}) one obtains
\begin{align}
2b_{13}^{(k-1)}\left( \tan\phi_{k} + \cot(2\phi_{k})\right) & = b_{11}^{(k)}-b_{22}^{(k)} \label{b13t}, \\
2b_{24}^{(k-1)}\left( \tan\psi_{k} + \cot(2\psi_{k})\right) & = b_{44}^{(k)}-b_{33}^{(k)}. \label{b24t}
\end{align}
Recall that for any $\zeta\in[-\frac{\pi}{4},\frac{\pi}{4}]\backslash\{0\}$ we have
$\displaystyle|\tan\zeta+ \cot2\zeta| = \frac{1}{2}\left| \tan\zeta+\frac{1}{\tan\zeta}\right|\geq1$. This implies
\begin{equation}\label{tancotan}
\min\big{\{}|\tan\phi_{k} + \cot(2\phi_{k})|,|\tan\psi_{k} + \cot(2\psi_{k})|\big{\}} \geq1, \quad k\geq 1.
\end{equation}
The relations (\ref{b13t}), (\ref{b24t}) and (\ref{tancotan}) imply
\begin{align}
2|b_{13}^{(k-1)}| \leq |b_{11}^{(k)}-b_{22}^{(k)}| & \leq |b_{11}^{(k)}-b_{44}^{(k)}|+|b_{22}^{(k)}-b_{44}^{(k)}|, \label{4x4a11a22a} \\
2|b_{24}^{(k-1)}| \leq |b_{44}^{(k)}-b_{33}^{(k)}| & \leq |b_{44}^{(k)}-b_{11}^{(k)}|+|b_{11}^{(k)}-b_{33}^{(k)}|. \label{4x4a11a22b}
\end{align}
After squaring and summing the inequalities (\ref{4x4a11a22a}) and (\ref{4x4a11a22b}), using (\ref{4x4delta}) and the inequality $(x+y)^{2}\leq1.5x^{2}+3y^{2}$
which holds for any real $x$ and $y$, for $k\geq1$ we get
\begin{equation}\label{intr3}
4\delta_{k-1}^{2} S^{2}(B) \leq3\left(|b_{11}^{(k)}-b_{44}^{(k)}|^{2}+|b_{22}^{(k)}-b_{44}^{(k)}|^{2}+|b_{11}^{(k)}-b_{33}^{(k)}|^{2} \right).
\end{equation}
Bounding the term $|b_{22}^{(k)}-b_{44}^{(k)}|^{2}+|b_{11}^{(k)}-b_{33}^{(k)}|^{2}$ is simple.
Using (\ref{4_ctg}), (\ref{4x4delta}) and Lemma~\ref{tm:lemma4}(v) we obtain
\begin{align}
|b_{22}^{(k)}-b_{44}^{(k)}|^{2}+|b_{11}^{(k)}-b_{33}^{(k)}|^{2} & = 2|b_{24}^{(k)}|^{2}|\cot2\psi_{k+1}|^{2}+2|b_{13}^{(k)}|^{2}|\cot2\phi_{k+1}|^{2} \nonumber \\
& \leq2\cdot4.49^{2}\delta_{k}^{2}\delta_{k+1}^{2}S^{2}(B). \label{4x4a22a44u}
\end{align}
The relations (\ref{intr3}) and (\ref{4x4a22a44u}) imply
\begin{equation}\label{dk-1}
\delta_{k-1}^{2} \leq\frac{3}{2}\cdot4.49^{2}\delta_{k}^{2}\delta_{k+1}^{2} + \frac{3}{4}\frac{|b_{11}^{(k)}-b_{44}^{(k)}|^{2}}{S^{2}(B)}, \quad k\geq1.
\end{equation}
Bounding $|b_{11}^{(k)}-b_{44}^{(k)}|^{2}$ is more demanding. From the relations (\ref{4_elementi}) for $k\geq1$ we have
\begin{align}
b_{11}^{(k+1)}-b_{33}^{(k+1)} & = b_{11}^{(k)}-b_{44}^{(k)} + b_{13}^{(k)}\tan\phi_{k+1}+b_{24}^{(k)}\tan\psi_{k+1} \nonumber \\
& = b_{11}^{(k)}-b_{44}^{(k)} + (b_{13}^{(k)}-b_{24}^{(k)})\tan\phi_{k+1} +b_{24}^{(k)} (\tan\phi_{k+1}+\tan\psi_{k+1}). \label{pom1}
\end{align}
From the relations (\ref{4_elementi}) we also get
\begin{equation}\label{pom2}
b_{13}^{(k)}-b_{24}^{(k)} = (b_{14}^{(k-1)}+b_{23}^{(k-1)})\cos(\phi_{k}+\psi_{k}), \quad k\geq1.
\end{equation}
Using (\ref{pom1}), (\ref{pom2}), (\ref{4_ctg}), Lemma~\ref{tm:lemma4}(iii), (\ref{4x4delta}), Lemma~\ref{tm:lemma4}(v), (\ref{nu}), (\ref{nu_k}) and (\ref{veliki_elementi_plus}) for $1\leq k\leq3$ one obtains
\begin{align}
|b_{11}^{(k)}-b_{44}^{(k)}| & \leq  |b_{11}^{(k+1)}-b_{33}^{(k+1)}|+|b_{13}^{(k)}-b_{24}^{(k)}|\cdot 1
+|b_{24}^{(k)}|\cdot|\tan\phi_{k+1}+\tan\psi_{k+1}| \nonumber\\
& \leq  2|b_{13}^{(k+1)}|\cdot|\cot2\phi_{k+2}| + |b_{14}^{(k-1)}+b_{23}^{(k-1)}|\cdot 1 + |b_{24}^{(k)}|4.444\delta_{k+1} \nonumber\\
& \leq 2\delta_{k+1}S(B)\cdot4.49\delta_{k+2} + \nu_{k-1} S(B) + \delta_kS(B)4.444\delta_{k+1} \nonumber\\
& \leq \left( 4.444\delta_k\delta_{k+1}+8.98\delta_{k+1}\delta_{k+2}+2.222^{k-1}\delta_{k-1}\cdots \delta_1\cdot \sqrt{2}\delta_1 \right) S(B).\label{pom3}
\end{align}

Here, for $k=1$ the term $2.222^{k-1}\delta_{k-1}\cdots\delta_{1}$ is replaced by one.
Next, we use the inequality $\sqrt{a^{2}+b^{2}}\leq|a|+|b|$ and combine (\ref{dk-1}) and (\ref{pom3}). After canceling by $S(B)$ we have
\begin{align*}
\delta_{k-1} & \leq \sqrt{3}\left[\frac{4.49}{\sqrt{2}}\delta_k\delta_{k+1}+\frac{1}{2}
\left(4.444\delta_k\delta_{k+1}+8.98\delta_{k+1}\delta_{k+2}+2.222^{k-1}\delta_{k-1}\cdots \delta_1\sqrt{2}\delta_1 \right)\right] \\
&\leq 9.348\delta_k\delta_{k+1}+7.777\delta_{k+1}\delta_{k+2}
+ 2.222^{k-1}\delta_{k-1}\cdots \delta_1\cdot \sqrt{1.5}\delta_1.
\end{align*}
Specifically, for $k=1,2,3$, one obtains
\begin{align*}
\delta_0 & \leq 1.225\delta_1+ 9.348\delta_1\delta_2+ 7.777\delta_{2}\delta_3, \\
\delta_1 & \leq 2.7214\delta_1^2+ 9.348\delta_2\delta_3+ 7.777\delta_{3}\delta_4, \\
\delta_2 & \leq 6.047\delta_2\delta_1^2 + 9.348\delta_3\delta_4+ 7.777\delta_{4}\delta_5.
\end{align*}
Since
\begin{align*}
(1-2.7214 \cdot\delta_1)^{-1} &< (1-2.7214\cdot 0.0045)^{-1} < 1.0124, \\
(1-6.047\cdot\delta_1^2)^{-1} &< (1-6.047\cdot 2\cdot 10^{-5} )^{-1} < 1.000121,
\end{align*}
it follows
$$\delta_1 < \delta_3 (9.464\delta_2+7.874\delta_4), \qquad \delta_2 < \delta_4 (9.35\delta_3+7.778\delta_5).$$
We used the fact that $\delta_{k}>0$ for $0\leq k\leq4$. Combining the bounds for $\delta_{0}$ and $\delta_{1}$ in order to eliminate the term $1.225\delta_{1}$ we obtain
$$\delta_0 <  9.348\delta_1\delta_2+19.371\delta_2\delta_3 + 9.646\delta_3\delta_4.$$
To bound $\delta_{2}$ by $\delta_{4}$ we just insert the upper bound $0.0045$ for $\delta_{3}$ and $\delta_{4}$ into the expression within the parentheses.
In a similar way (as $(9.35+7.778)\cdot2\cdot10^{-5}$) the bound for $\delta_{2}$ can be obtained.
To bound $\delta_{1}$ by $\delta_{3}$ we use the bounds $0.000347$ and $0.0045$ for $\delta_{2}$ and~$\delta_{4}$.
Finally, the upper bounds for $\delta_{1}$ and $\delta_{0}$ are obtained (in this order) by inserting the best available bounds into the appropriate expressions.
\qed

\subsection{Proof of Lemma~\ref{tm:lemma7}}

\begin{itemize}
\item[(i)] The assertion follows from (\ref{pom2}), (\ref{nu}) and (\ref{nu_k}) or (\ref{veliki_elementi_plus}).
\item[(ii)] We use the parallelogram law, $(u+v)^{2}+(u-v)^{2}=2(u^{2}+v^{2})$ for $u,v$ real, and the definition of $\delta_{k}$ from (\ref{4x4delta}).
For $k=3,4$ the inequalities (\ref{a13+a24-3}) follow from (\ref{diffpivel}) and Lemma~\ref{tm:lemma5}. We have
\begin{align*}
\frac{|b_{13}^{(3)}+b_{24}^{(3)}|}{S(B)} &> \sqrt{2\delta_3^2 - 6.9824^2\delta_2^2\delta_1^4} > \sqrt{2 - 6.9824^2\delta_2^2\delta_1^2\cdot 0.0388^2}\delta_3 \\
&> \sqrt{2 - 0.0734\cdot 34.7^2\cdot 17.5^2\cdot 10^{-20}}\delta_3 \geq 1.41421\delta_3, \\
\frac{|b_{13}^{(4)}+b_{24}^{(4)}|}{S(B)} &> \sqrt{2\delta_4^2 - 15.52^2\delta_3^2\delta_2^2\delta_1^4} >
\sqrt{2- 15.52^2\cdot 0.0771^2\delta_3^2\delta_1^4}\delta_4\\
&> \sqrt{2 - 1.432\cdot 0.0045^2\cdot 17.5^4\cdot 10^{-20}}\delta_4 > 1.41421\delta_4.
\end{align*}
\item[(iii)] We use the first case in Lemma~\ref{tm:lemma4}(i)
\begin{align}
|\cot2\phi_{k} - \cot2\psi_{k}| & = \frac{|\sin 2(\phi_k -\psi_k)|}{|\sin 2\phi_k \sin 2\psi_k|} \geq |\sin 2(\phi_k -\psi_k)| \nonumber \\
& =  2|\cos (\alpha_k +\beta_k)|\cos (\phi_k -\psi_k) > \left(2-4.937284\delta_k^2\right)\cos (\phi_k -\psi_k) \nonumber \\
&> 1.9999\cos(\phi_k -\psi_k), \quad 1\leq k\leq 5,\label{c-}
\end{align}
with $\delta_{k}^{2} < 2\epsilon=2\cdot10^{-5}$. The proof of (\ref{c-}) is the same if we use $\alpha_{k}'$, $\beta_{k}'$ instead of $\alpha_{k}$, $\beta_{k}$, respectively.
From the relations (\ref{4_elementi}) it follows that
\begin{equation}\label{pom2_b}
b_{13}^{(k)}+b_{24}^{(k)} = (b_{14}^{(k-1)}-b_{23}^{(k-1)})\cos(\phi_{k}-\psi_{k}), \quad k\geq1.
\end{equation}
Hence, for $3\leq k\leq4$ we can use the assertion (ii) to obtain
\begin{align*}
\cos(\phi_k -\psi_k) & = \frac{|b_{13}^{(k)}+b_{24}^{(k)}|}{|b_{14}^{(k-1)}-b_{23}^{(k-1)}|}
\geq \frac{1.41421\delta_k S(B)}{\sqrt{2}\sqrt{(b_{14}^{(k-1)})^2+(b_{23}^{(k-1)})^2}} \\
& > \frac{1.41421\delta_k S(B)}{\sqrt{2}S(B)} > 0.999997\delta_k.
\end{align*}
Combining that with (\ref{c-}) one obtains $|\cot2\phi_{k} - \cot2\psi_{k}| > 1.999894\delta_{k}$, $3\leq k\leq4$.
\qed
\end{itemize}

\subsection{Proof of Lemma~\ref{tm:lemma5} in the case $|b_{14}-b_{23}|\leq\sqrt{2}\delta_{1}S(B)$}

The relations (\ref{4_ctg})--(\ref{dk-1}) remain the same except for the relation (\ref{veliki_elementi_plus}) in which $\nu_{k}$ and $\nu$ are replaced by $\nu_{k}^{-}$ and $\nu^{-}=\nu_{0}^{-}$, respectively. For $k\geq1$ the relation (\ref{pom1}) can be written as
\begin{equation}\label{pom1_b}
b_{11}^{(k+1)}-b_{33}^{(k+1)} = b_{11}^{(k)}-b_{44}^{(k)} + (b_{13}^{(k)}+b_{24}^{(k)})\tan\phi_{k+1} +b_{24}^{(k)} (\tan\psi_{k+1}-\tan\phi_{k+1}),
\end{equation}
and instead of (\ref{pom2}) we use (\ref{pom2_b}).
The relation (\ref{pom3}) takes the form
\begin{align*}
|b_{11}^{(k)}-b_{44}^{(k)}| & \leq  |b_{11}^{(k+1)}-b_{33}^{(k+1)}|+|b_{13}^{(k)}+b_{24}^{(k)}|\cdot1 + |b_{24}^{(k)}|\cdot|\tan\psi_{k+1}-\tan\phi_{k+1}| \\
&\leq  2|b_{13}^{(k+1)}|\cdot|\cot2\phi_{k+2}| + |b_{14}^{(k-1)}-b_{23}^{(k-1)}|\cdot 1 + |b_{24}^{(k)}|4.444\delta_{k+1} \\
&\leq 2\delta_{k+1}S(B)\cdot4.49\delta_{k+2} + \nu_{k-1}^{-} S(B) + \delta_kS(B)4.444\delta_{k+1} \\
&\leq \left(4.444\delta_k\delta_{k+1}+8.98\delta_{k+1}\delta_{k+2}+2.222^{k-1}\delta_{k-1}\cdots \delta_1\cdot \sqrt{2}\delta_1\right) S(B),
\end{align*}
for $1\leq k\leq3$. Here we used the relations (\ref{pom1_b}), (\ref{pom2_b}), (\ref{nu_k-}) and the assertions (iii) and (v) of Lemma~\ref{tm:lemma4_b}.
The rest of the proof follows the remaining lines in the proof of Lemma~\ref{tm:lemma5}.
\qed

\subsection{Proof of Lemma~\ref{tm:lemma7_b}}

The proof of the first two assertions is quite similar to the proof of the appropriate assertions of Lemma~\ref{tm:lemma7}.
To prove the third assertion instead of the relation (\ref{c-}) we now have
\begin{align}
|\cot2\phi_{k} + \cot2\psi_{k}| & = \frac{|\sin 2(\phi_k +\psi_k)|}{|\sin 2\phi_k \sin 2\psi_k|} \geq |\sin 2(\phi_k +\psi_k)| \nonumber \\
& = 2|\cos (\alpha_k +\beta_k)|\cos(\phi_k +\psi_k) > \left(2-4.937284\delta_k^2\right)\cos(\phi_k +\psi_k) \nonumber \\
& > 1.9999\cos(\phi_k +\psi_k), \quad 1\leq k\leq5, \label{c--}
\end{align}
where we used $\delta_{k}^{2} < 2\epsilon=2\cdot10^{-5}$. The proof of (\ref{c-}) is the same if we use $\alpha_{k}'$, $\beta_{k}'$ instead of $\alpha_{k}$, $\beta_{k}$, respectively. Using (\ref{pom2}) and the assertion (ii) for $3\leq k\leq4$ we obtain
\begin{align*}
\cos(\phi_k +\psi_k) & = \frac{|b_{13}^{(k)}-b_{24}^{(k)}|}{|b_{14}^{(k-1)}+b_{23}^{(k-1)}|}
\geq \frac{1.41421\delta_kS(B)}{\sqrt{2}\sqrt{(b_{14}^{(k-1)})^2+(b_{23}^{(k-1)})^2}} 
 > \frac{1.41421\delta_kS(B)}{\sqrt{2}S(B)} > 0.999997\delta_k.
\end{align*}
Combining that with (\ref{c--}) we get
$|\cot2\phi_{k} + \cot2\psi_{k}| > 1.999894\delta_{k}$, $3\leq k\leq4$.
\qed

\section{Proofs related to Section~\ref{sec5}}\label{appC}

\subsection{Proof of Proposition~\ref{tm:prop5.1}}

We use the operator $\mathcal{T}$ from Definition~\ref{T}. Let $B=B^{0}=H(\epsilon)$ and $B^{(k)}=\mathcal{T}^{k}(H)$ for $k\geq1$.

For $k=1$ we compute $B^{(1)}$ from $B$. The elements $b_{13}$ and $b_{24}$ are annihilated and the off-norm reduction equals
\begin{equation}\label{prr1}
S^{2}(B)-S^{2}(B^{(1)})=(2\epsilon)^{2}+\epsilon^{2}=5\epsilon^{2}.
\end{equation}

For the rotation angels $\phi_{1}$ and $\psi_{1}$ we have
\begin{align*}
\tan(2\phi_{1}) & =\frac{4\epsilon}{\epsilon^{1.5}}=\frac{4}{\sqrt{\epsilon}}, \quad 0<\phi_{1}<\frac{\pi}{4}, \\
\tan(2\psi_{1}) & =-\frac{2\epsilon}{\epsilon^{1.5}}=-\frac{2}{\sqrt{\epsilon}}, \quad -\frac{\pi}{4}<\psi_{1}<0.
\end{align*}
Using the notation from Lemma~\ref{tm:lemma4} we have
\begin{equation}\label{prkutevi1}
\phi_{1}=\frac{\pi}{4}-\alpha_{1}, \qquad \psi_{1}=-\frac{\pi}{4}+\beta_{1}, \qquad \alpha_{1}>0, \ \beta_{1}>0.
\end{equation}
Hence
\begin{align}
\tan(2\alpha_1) & = \frac{1}{\cot(2\alpha_1)} = \frac{1}{\tan(\frac{\pi}{2}-2\alpha_1)}=\frac{1}{\tan (2\phi_1)}=\frac{\sqrt{\epsilon}}{4}, \label{pr:_t2a} \\
\tan(2\beta_1) & = \frac{1}{\cot(2\beta_1)} = \frac{1}{\tan(\frac{\pi}{2}-2\beta_1)}=\frac{1}{-\tan(2\psi_1)}=\frac{\sqrt{\epsilon}}{2}, \label{pr:_t2b}
\end{align}
and
\begin{equation}\label{prtilde1}
\alpha_{1}=\frac{1}{2}\arctan(\frac{\sqrt{\epsilon}}{4}), \qquad \beta_{1}=\frac{1}{2}\arctan(\frac{\sqrt{\epsilon}}{2}).
\end{equation}
Since $\displaystyle\arctan(x)=x-\frac{x^{3}}{3}+\frac{x^{5}}{5}-\frac{x^{7}}{7}+\cdots$, for $0<x\leq1$ we have
$$x-\frac{x^3}{3} < \arctan (x) < x-\frac{x^3}{3}+\frac{x^5}{5}<x.$$
The relation (\ref{prtilde1}) implies
\begin{align*}
0<\frac{\sqrt{\epsilon}}{8}-\frac{\epsilon\sqrt{\epsilon}}{384} & < \alpha_1 < \frac{\sqrt{\epsilon}}{8}-\frac{\epsilon\sqrt{\epsilon}}{384} + \frac{\epsilon^2\sqrt{\epsilon}}{10240}
<\frac{\sqrt{\epsilon}}{8}, \\
0<\frac{\sqrt{\epsilon}}{4}-\frac{\epsilon\sqrt{\epsilon}}{48} & < \beta_1 <
\frac{\sqrt{\epsilon}}{4}-\frac{\epsilon\sqrt{\epsilon}}{48}+\frac{\epsilon^2\sqrt{\epsilon}}{320}< \frac{\sqrt{\epsilon}}{4},
\end{align*}
and consequently
\begin{equation}\label{pr1tmv}
\frac{3}{8}\sqrt{\epsilon}\left(1-\frac{\epsilon}{16}\right) < \alpha_1+\beta_1
< \frac{3}{8}\sqrt{\epsilon}- \frac{3}{128}\epsilon\sqrt{\epsilon}\left(1-\frac{11}{80}\epsilon \right)\epsilon\sqrt{\epsilon}< \frac{3}{8}\sqrt{\epsilon}.
\end{equation}
Since $\displaystyle\phi_{1}-\psi_{1} = \frac{\pi}{2}-(\alpha_{1}+\beta_{1})$, we have $\cos(\phi_{1}-\psi_{1})=\sin(\alpha_{1}+\beta_{1})$.
Therefore, using the relation (\ref{pr1tmv}) we obtain
\begin{align*}
\sin(\alpha_{1}+\beta_{1}) & < \alpha_{1}+\beta_{1} < \frac{3}{8}\sqrt{\epsilon}, \\
\sin(\alpha_{1}+\beta_{1}) & > (\alpha_{1}+\beta_{1})\left(1-\frac{(\alpha_{1}+\beta_{1})^{2}}{6}\right) >
\frac{3}{8}\sqrt{\epsilon}\left( 1-\frac{\epsilon}{16}\right) \left(1-\frac{9\epsilon}{6\cdot64}\right) > \frac{3}{8}\sqrt{\epsilon} - \frac{33}{1024}\epsilon\sqrt{\epsilon}.
\end{align*}
Thus,
\begin{equation}\label{prcos1}
\frac{3}{8}\sqrt{\epsilon}\left( 1-\frac{11}{128}\epsilon\right) < \cos(\phi_{1}-\psi_{1})< \frac{3}{8}\sqrt{\epsilon}.
\end{equation}
Note that
\begin{equation}\label{prr2}
S^{2}(B^{(1)})-S^{2}(B^{(2)}) = (b_{13}^{(1)})^{2}+(b_{24}^{(1)})^{2},
\end{equation}
hence we have to bound $(b_{13}^{(1)})^{2}+(b_{24}^{(1)})^{2}$.
From the relations (\ref{4_elementi}) it follows that
\begin{align}
b_{13}^{(1)} & = (\epsilon-1)\cos{\phi_1}\cos{\psi_1}-\sin{\phi_1}\sin{\psi_1} = \epsilon\cos\phi_1\cos\psi_1-\cos(\phi_1-\psi_1), \label{prb13_1} \\
b_{24}^{(1)} & = (\epsilon-1)\sin{\phi_1}\sin{\psi_1} -\cos{\phi_1}\cos{\psi_1} = \epsilon\sin\phi_1\sin\psi_1 -\cos(\phi_1 -\psi_1). \label{prb24_1}
\end{align}
Inserting (\ref{prb13_1}) and (\ref{prb24_1}) into (\ref{prr2}) and using (\ref{prcos1}) we obtain
\begin{align}
\nonumber (b_{13}^{(1)})^2+(b_{24}^{(1)})^2 & = (\epsilon\cos\phi_1\cos\psi_1-\cos(\phi_1-\psi_1))^2 + (\epsilon\sin\phi_1\sin\psi_1-\cos(\phi_1-\psi_1))^2 \\
\nonumber & = 2\cos^2(\phi_1-\psi_1) - 2\epsilon\cos(\phi_1-\psi_1)(\cos\phi_1\cos\psi_1+\sin\phi_1\sin\psi_1) \\
\nonumber & \qquad + \epsilon^2 (\cos^2\phi_1\cos^2\psi_1+\sin^2\phi_1\sin^2\psi_1) \leq  2\cos^2(\phi_1-\psi_1)(1-\epsilon ) + \epsilon^2 \\
\label{prr22} & \leq  \frac{9}{31}(1+\frac{22}{9}\epsilon) \epsilon <0.29033\epsilon <\frac{3}{10}\epsilon.
\end{align}
Let $k=2$. We have
\begin{equation}\label{prtan2}
\tan(2\phi_{2})=\frac{2b_{13}^{(1)}}{b_{11}^{(1)}-b_{33}^{(1)}}, \qquad \tan(2\psi_{2})=\frac{2b_{24}^{(1)}}{b_{22}^{(1)}-b_{44}^{(1)}}.
\end{equation}
From the relations (\ref{prb13_1}), (\ref{prb24_1}) and (\ref{prcos1}) we conclude that the pivot elements $b_{13}^{(1)}$ and $b_{24}^{(1)}$ are negative.
Using the relations (\ref{prkutevi1}) and (\ref{prcos1}) we bound their moduli from below
\begin{align*}
|b_{13}^{(1)}| & = -b_{13}^{(1)} = \cos(\phi_1-\psi_1) - \epsilon\cos\phi_1\cos\psi_1 \\
& > \frac{3}{8}\sqrt{\epsilon}\left(1-\frac{11}{128}\epsilon\right) - \frac{\epsilon}{2}(\cos(\alpha_1)+\sin (\alpha_1))(\cos(\beta_1)+\sin (\beta_1)) \\
& > \frac{3}{8}\sqrt{\epsilon} - \frac{33}{1024}\epsilon\sqrt{\epsilon}-\frac{\epsilon}{2}(1+\frac{\sqrt{\epsilon}}{8})
(1+\frac{\sqrt{\epsilon}}{4}) > \frac{3}{8}\sqrt{\epsilon}- 0.5007\epsilon,\\
|b_{24}^{(1)}| & = -b_{24}^{(1)} = \cos(\phi_1-\psi_1) - \epsilon\sin\phi_1\sin\psi_1 \\
& \geq \frac{3}{8}\sqrt{\epsilon} - \frac{33}{1024}\epsilon\sqrt{\epsilon} +\frac{\epsilon}{2}(\cos(\alpha_1)-\sin (\alpha_1))(\cos(\beta_1)-\sin (\beta_1)) > \frac{3}{8}\sqrt{\epsilon}.
\end{align*}
Moreover, using (\ref{4_elementi}) or (\ref{pom1}) and (\ref{pr:_t2a}), (\ref{pr:_t2b}) we obtain
\begin{align*}
b_{11}^{(1)}-b_{33}^{(1)} & = b_{11}-b_{44} + b_{13}\tan{\phi_1}+b_{24}\tan{\psi_1} = \epsilon+\epsilon^{1.5}+2\epsilon\tan{\phi_1}-\epsilon\tan{\psi_1} \\
& = \epsilon + \epsilon\sqrt{\epsilon} + 2\epsilon\tan (\frac{\pi}{4}-\alpha_1) + \epsilon\tan (\frac{\pi}{4}-\beta_1) = \epsilon + \epsilon\sqrt{\epsilon}+ 2\epsilon\frac{1-\tan (\alpha_1)}{1+\tan (\alpha_1)}+\epsilon\frac{1-\tan (\beta_1)}{1+\tan (\beta_1)} \\
&\geq \epsilon + \epsilon\sqrt{\epsilon}+ 2\epsilon \left(1-2\tan (\alpha_1)\right) +\epsilon \left(1-2\tan (\beta_1)\right) \geq 4\epsilon + \epsilon\sqrt{\epsilon}
-4\epsilon \tan (\alpha_1)-2\epsilon\tan (\beta_1) \\
& \geq 4\epsilon + \epsilon\sqrt{\epsilon} -4\epsilon\cdot \frac{1}{2}\frac{\sqrt{\epsilon}}{4} -2\epsilon\cdot\frac{1}{2}\frac{\sqrt{\epsilon}}{2}=4\epsilon > 0,
\end{align*}
and
\begin{align*}
b_{22}^{(1)}-b_{44}^{(1)} & = b_{33}-b_{22}-(b_{13}\tan{\phi_1}+b_{24}\tan{\psi_1}) = \epsilon-\epsilon\sqrt{\epsilon} -(2\epsilon\tan{\phi_1}-\epsilon\tan{\psi_1}) \\
& \leq \epsilon-\epsilon\sqrt{\epsilon}-\big{(}2\epsilon(1-2\tan(\alpha_1)) + \epsilon(1-2\tan(\beta_1))\big{)} \\
& = -2\epsilon-\epsilon\sqrt{\epsilon}+ 4\epsilon \tan (\alpha_1)+2\epsilon\tan (\beta_1) \\
& \leq -2\epsilon-\epsilon\sqrt{\epsilon}+ 4\epsilon\cdot \frac{1}{2}\frac{\sqrt{\epsilon}}{4} +2\epsilon\cdot\frac{1}{2}\frac{\sqrt{\epsilon}}{2}
=-2\epsilon < 0.
\end{align*}
Hence, we conclude that $\phi_{2}<0$, $\psi_{2}>0$. Like in Lemma \ref{tm:lemma4}(i), we set
$$\phi_2=-\frac{\pi}{4}+\alpha_2', \quad \psi_2=\frac{\pi}{4}-\beta_2', \quad \alpha_2'>0, \ \beta_2'>0.$$
From the relations (\ref{prtan2}), (\ref{pr:_t2a}) and (\ref{pr:_t2b}) we obtain
\begin{align*}
|\cot(2\phi_2)| & = \frac{b_{11}^{(1)}-b_{33}^{(1)}}{2|b_{13}^{(1)}|} \leq \frac{\epsilon+\epsilon\sqrt{\epsilon}+2\epsilon (1-\tan (\alpha_1)) +\epsilon (1-\tan (\beta_1))}{\frac{3}{4}\sqrt{\epsilon}-1.0014\epsilon} \\
& \leq  \frac{4\epsilon+\epsilon\sqrt{\epsilon}}{(0.75-1.0014\cdot 0.0032)\sqrt{\epsilon}} < \frac{4 +0.0032}{0.74679552}\sqrt{\epsilon} = 5.361\sqrt\epsilon, \\
|\cot(2\psi_2)| & =\frac{|b_{22}^{(1)}-b_{44}^{(1)}|}{2|b_{24}^{(1)}|}
=\frac{2\epsilon\tan{\phi_1}-\epsilon\tan{\psi_1}-(\epsilon-\epsilon\sqrt{\epsilon})}{\frac{3}{4}\sqrt{\epsilon}} \\
& \leq \frac{2\epsilon+\epsilon\sqrt{\epsilon}}{0.75\sqrt{\epsilon}} < \frac{2 + 0.0032} {0.75}\sqrt{\epsilon} < 2.671\sqrt\epsilon.
\end{align*}
Then, for $\alpha_{2}'$ and $\beta_{2}'$ it holds
\begin{align*}
\tan(2\alpha_2') & =\tan(\frac{\pi}{2}+2\phi_2)=|\cot(2\phi_2)|< 5.361\sqrt\epsilon, \\
\alpha_2' & <\frac{1}{2}\arctan(5.361\sqrt{\epsilon})< 2.681\sqrt{\epsilon}, \\
\tan(2\beta_2') & =\tan(\frac{\pi}{2}-2\psi_2)=\cot(2\psi_2)< 2.671\sqrt\epsilon, \\
\beta_2' & <\frac{1}{2}\arctan(2.671\sqrt{\epsilon})< 1.336\sqrt{\epsilon}.
\end{align*}

Next we bound the off-norm reduction in the third parallel step which equals $(b_{13}^{(2)})^{2}+(b_{24}^{(2)})^{2}$.
We use the relations (\ref{prva}), (\ref{new1}) and (\ref{b14pmb23}) to obtain
\begin{align}
(b_{13}^{(2)})^2+(b_{24}^{(2)})^2 & = \big{(}(b_{14}^{(1)})^2+(b_{23}^{(1)})^2\big{)}\cos^2(\phi_2-\psi_2) - 2(b_{14}^{(1)}+b_{23}^{(1)})^2\cdot \frac{1}{4}\sin{2\phi_2}\sin{2\psi_2}   \nonumber \\
& = S^2(B^{(2)})\cos^2(\phi_2-\psi_2) + \frac{1}{2}\sin^2 (\phi_1+\psi_1)|\sin{2\phi_2}\sin{2\psi_2}|(b_{14}+b_{23})^2 \nonumber \\
& \leq S^2(B)\cos^2(\alpha_2'+\beta_2'-\frac{\pi}{2}) + \frac{1}{2}\cdot \sin^2 (\beta_1-\alpha_1) \cdot 1\cdot\epsilon^2 \nonumber \\
& \leq 2(\alpha_2'+\beta_2')^2+\frac{1}{32}\epsilon^3 \leq \left(2\cdot 4.017^2 +\frac{10^{-10}}{32}\right)\epsilon \leq 32.273\epsilon. \label{prr3}
\end{align}
Finally, from (\ref{prr1}), (\ref{prr2}), (\ref{prr22}) and (\ref{prr3}) it follows
\begin{align*}
S^2(H^{(0)})-S^2(H^{(6)}) & = S^2(B)-S^2(B^{(3)}) = b_{13}^2+b_{24}^2 + (b_{13}^{(1)})^2+(b_{24}^{(1)})^2 + (b_{13}^{(2)})^2+(b_{24}^{(2)})^2 \\
& < 5\epsilon^2+0.29033\epsilon+32.273\epsilon < 32.56338\epsilon < \frac{32.56338}{2-2\epsilon}\epsilon S^2(H^{(0)}) \\
& < 16.282\epsilon S^2(H^{(0)}).
\qed
\end{align*}

\end{document}